\documentclass[12pt]{amsart}
\usepackage{enumerate}
\usepackage{xcolor}
\usepackage[normalem]{ulem}
\usepackage{letltxmacro}
\usepackage[all]{xy}
\usepackage{comment}
\usepackage{mathalpha}
\usepackage{mathrsfs}
\usepackage{tikz-cd}
\usepackage{amsmath,amssymb,setspace,nicefrac, yhmath, amscd}
\usepackage[active]{srcltx}
\usepackage[pagebackref, colorlinks, linkcolor=blue, citecolor=blue, urlcolor=blue, hypertexnames=true]{hyperref}
\usepackage{amsrefs}
\usepackage{dirtytalk}
\usepackage{xcolor}
\allowdisplaybreaks
\usepackage[all]{xy}
\newcommand{\setbuilder}[2] { \left\{ #1 \enskip \middle| \enskip #2 \right\} }

\let\oldaddcontentsline\addcontentsline
\newcommand{\stoptocentries}{\renewcommand{\addcontentsline}[3]{}}
\newcommand{\starttocentries}{\let\addcontentsline\oldaddcontentsline}
\theoremstyle{plain}
\newtheorem{theorem}{Theorem}[section]
\newtheorem{corollary}[theorem]{Corollary}
\newtheorem{proposition}[theorem]{Proposition}
\newtheorem{lemma}[theorem]{Lemma}
\theoremstyle{definition}
\newtheorem{remark}[theorem]{Remark}
\newtheorem{example}[theorem]{Example}

\newtheorem{definition}[theorem]{Definition}

\newtheorem{claim}{Claim}

\newtheorem{maintheorem}{Theorem}

\theoremstyle{plain}

\newcommand{\bE}{{\mathbb{E}}}

\def\Fix{\mr{Fix}}
  \newcommand{\N}{{\mathcal{N}}}

\def\Ypi{(Y,\pi_Y)}
\def\mr{\mathrm}
\def\mb{\mathbb}

\def\pF{\partial_F}
\def\Inn{\mr{Inn}}
\newcommand{\Prob}{\mathcal{P}} % Space of probability measures
\newcommand{\al}{\alpha}

\renewcommand{\phi}{\varphi}
\newcommand{\upchi}{{\raise.35ex\hbox{\ensuremath{\chi}}}}

\def\pY{\pi_Y}
\newcommand{\Ad}{\operatorname{Ad}}

\newcommand{\Aut}{\operatorname{Aut}}

\newcommand{\id}{{\operatorname{id}}}

\def\Prob{\mathrm{Prob}}

\newcommand{\supp}{\operatorname{supp}}

\usepackage{xcolor}

\usepackage{hyperref}
\title[Rigidity of Generalized Furstenberg Boundaries]{Rigidity of Generalized Furstenberg Boundaries and Applications to Intermediate Crossed Products}
\author[Amrutam]{Tattwamasi Amrutam}
\address{Institute of Mathematics of the Polish Academy of Sciences, ul. Sniadeckich 8, 00-656, Warszawa, Poland}
\email{tattwamasiamrutam@gmail.com}
\author[Liu]{Chunlin Liu}\address{School of Mathematical Sciences, Dalian University of Technology, Dalian, 116024, P.R. China, and Institute of Mathematics of the Polish Academy of Sciences, ul. Sniadeckich 8, 00-656, Warszawa, Poland}
\email{chunlinliu@mail.ustc.edu.cn}
\date{\today}
\begin{document}
\begin{abstract}
We develop a relative boundary theory for actions of discrete groups on
compact spaces and use it to derive rigidity results for reduced crossed
products. For a discrete group $\Gamma$ acting on a compact space $X$ and a
subgroup $H$, we construct a universal boundary over $X$ which is minimal as
a $\Gamma$-system and strongly proximal with respect to $H$. When $H\le_c\Gamma$ is commensurated
and the $H$-action on $X$ is minimal, we show that this universal boundary
agrees, in a canonical $\Gamma$-equivariant way, with the generalized
Furstenberg boundary of $(H,X)$, thereby unifying and extending earlier
results on relative boundaries.

As an application, we introduce the notion of an $X$-plump subgroup given a $\Gamma$-space $X$,
a generalized version of plumpness tailored to crossed products. Under natural
dynamical hypotheses, this leads to new examples of irreducible $C^*$-inclusions. Under additional assumptions, we also show that every intermediate
$C^*$-algebra is a crossed product.
\end{abstract}
\dedicatory{Dedicated to Prof. Swadheen Pattanayak on his 85th birthday with great admiration.}

\maketitle
\tableofcontents
\section{Introduction}
Boundaries provide a unifying framework for studying the interplay between group actions, rigidity phenomena, and operator-algebraic structure.
For a discrete group $\Gamma$, the Furstenberg boundary $\partial_F\Gamma$ is characterized as the universal minimal strongly proximal $\Gamma$-space, and it has become a fundamental object in the analysis of $C^*$-simplicity~\cite{kalantar2017boundaries}, ideal structure of crossed products~\cite{breuillard2017c}, and boundary amenability.
A growing body of work has demonstrated that subtle algebraic properties of $C_r^*(\Gamma)$ are reflected in the dynamics of its action on $\partial_F\Gamma$ and its topological factors~\cite{ozawa2025proximality}.

In many situations of interest, one is not given a group in isolation, but rather a group action $\Gamma\curvearrowright X$ together with a subgroup $H\leq \Gamma$ encoding additional geometric or algebraic structure.
Examples include Hecke pairs ($(G,H)$ is called a Hecke pair if every double
coset of $H$ is a union of finitely many left cosets of $H$), lattices in locally compact groups, stabilizers of points or subspaces, and actions arising from arithmetic or geometric constructions.
In such settings, it is natural to ask how boundary phenomena behave \emph{relative} to the base space $X$ and how they interact with the subgroup structure of $\Gamma$.

Following ideas of Glasner, the \emph{generalized Furstenberg boundary} $\partial_F(\Gamma,X)$ is defined as the universal minimal strongly proximal extension of the $\Gamma$-space $X$~\cite{SGlasner1975,naghavi2020furstenberg}.
This construction recovers the classical Furstenberg boundary when $X$ is a point, but in general it encodes the strongest possible proximal behavior compatible with the dynamics on $X$.
Generalized Furstenberg boundaries have appeared in several contexts, particularly in recent works of Naghavi~\cite{naghavi2020furstenberg} and Kawabe~\cite{kawabe2017uniformlyrecurrentsubgroupsideal}, where freeness of the action on $\partial_F(\Gamma,X)$ was shown to characterize simplicity of the reduced crossed product $C(X)\rtimes_r \Gamma$ under minimality of the action $\Gamma\curvearrowright X$.
Despite these developments, a systematic theory describing how generalized boundaries behave when transitioning between groups and subgroups—analogous to the Furstenberg boundary—has remained undeveloped. As such, our first main result is the existence of a universal generalized $H$-boundary for any subgroup $H\le\Gamma$ (see Section~\ref{sec:pre} for definitions).
\begin{maintheorem}\label{thmA}
	Let $H\le \Gamma$. Let  $X$ be a $\Gamma$-space which is $H$-minimal.
	Consider the class $\mathcal C$ of all $\Gamma$-extensions $\pi:Y\to X$ such that
	\begin{enumerate}
		\item $Y$ is $\Gamma$-minimal;
		\item $\pi$ is an $H$-strongly proximal extension over $X$.
	\end{enumerate}
	Morphisms in $\mathcal C$ are $\Gamma$-equivariant continuous maps $f:Y\to Y'$
	satisfying $\pi' \circ f=\pi$.
	
Then there exists a universal object $\pi_{RB}: RB_X(\Gamma,H)\to X$ in $\mathcal C$:
for every $\pi:Y\to X$ in $\mathcal C$, there exists a $\Gamma$-morphism
$f: RB_X(\Gamma,H)\to Y$ such that $\pi\circ f=\pi_{RB}$.
\end{maintheorem}
Now, given $H\le \Gamma$, can one identify $RB_X(\Gamma,H)$ with the generalized
Furstenberg boundary $\partial_F(H,X)$ (as $\Gamma$--extensions over $X$)?
  To answer this question, we must first ask if the generalized boundary $\partial_F(H,X)$ can be promoted to a generalized $\Gamma$-boundary?
When $X$ is a one-point space, Li and Scarparo~\cite{li2023c} (also see~\cite{dai2019universal}) showed that the Furstenberg boundary of $\Gamma$ can be reconstructed from that of $H$ via a canonical extension of the action, when $H\le_c\Gamma$ is commensurable. Related techniques have also played a key role in boundary-based approaches to $C^*$-simplicity due to Breuillard, Kalantar, Kennedy, and Ozawa~\cite{breuillard2017c}, where such an extension was first shown for a normal subgroup $N\triangleleft\Gamma$ with respect to the Furstenberg boundary.

In this paper, we generalize the result of Li and Scarparo to that of a generalized Furstenberg boundary. Recall that $H\leq \Gamma$ is called a \emph{commensurated subgroup}, if  $H\cap gHg^{-1}$ has finite index in both $H$ and $gHg^{-1}$ for all $g\in\Gamma$.
While commensurated subgroups share certain large-scale features with normal subgroups, they lack equivariance of conjugation, and as a consequence classical extension arguments do not apply.
The main technical difficulty is to construct a $\Gamma$-action on $\partial_F(H,X)$ that is compatible with the given $H$-action and is independent of auxiliary finite-index choices.

Our second main result shows that, under a natural minimality assumption on the base action, this difficulty can be overcome.
We prove that the generalized Furstenberg boundary $\partial_F(H,X)$ admits a canonical and unique extension of the $H$-action to $\Gamma$, and that with this extension, $\Gamma\curvearrowright\partial_F(H,X)$ is a relative boundary over $X$, thus identifying it with $RB_X(\Gamma,H)$.
\begin{maintheorem}\label{thmB}
Let  $X$ be a minimal $\Gamma$-space, and let $H\leq \Gamma$ be a commensurated subgroup.
Then the generalized Furstenberg boundary $\partial_F(H,X)$ admits a unique extension of the $H$-action to a minimal $\Gamma$-action.
Moreover, with this extended action, $\partial_F(H,X)$ is $H$-equivariantly homeomorphic to $\partial_F(\Gamma,X)$.
\end{maintheorem}
The proof relies on rigidity phenomena for minimal strongly proximal actions under finite-index restrictions, together with a careful analysis of canonical extensions along intersections of conjugates of $H$.   

When $X$ is an $H$-minimal $\Gamma$-space and $H\le_c\Gamma$
is commensurated, the universal
generalized $H$-boundary $RB_X(\Gamma,H)$, constructed in
Theorem~\ref{thmA} agrees with the generalized Furstenberg boundary
$\partial_F(H,X)$ endowed with the extended $\Gamma$-action from
Theorem~\ref{thmB} (see Theorem~\ref{thm:same}). In particular, there is a $\Gamma$-equivariant
homeomorphism $RB_X(\Gamma,H)\cong \partial_F(H,X)$.

In the second part of the paper, we investigate operator-algebraic consequences of the boundary extension theory.
Our approach is motivated by the Powers averaging principle, which underlies the simplicity of the reduced group $C^*$-algebra of discrete groups and has since become a central mechanism in the study of $C^*$-simplicity~\cite{haagerup, Ken}. Building on this idea, the first-named author introduced the notion of \emph{plump subgroups}~\cite{amrutam2021intermediate}, providing a dynamical averaging criterion given by a subgroup of the ambient group. More recently, the first-named author and Ursu ~\cite{amrutam2022generalized} developed a generalized Powers averaging framework, extending this mechanism to crossed products of commutative $C^*$-algebras.

Guided by these developments, we introduce the notion of an \emph{$X$-plump} subgroup (see Definition~\ref{def:genplump}), which captures a form of relative averaging compatible with a given $\Gamma$-space $X$ and its generalized Furstenberg boundary.
This notion can be viewed as a relative analogue of plumpness, adapted to crossed products of the form $C(X)\rtimes_r \Gamma$ (see Section~\ref{subsec:crossed} for more on the reduced crossed products). 
Using $X$-plumpness and the relationship between generalized Furstenberg boundaries of an ambient group and those of its subgroup, we obtain new examples of irreducible $C^*$-inclusions.
\begin{maintheorem}\label{thmC}
    Let $\pi: X \to Y$ be a continuous $\Gamma$-equivariant onto map, where $X,Y$ are $\Gamma$-spaces. Suppose that there exists a subgroup $H\le \Gamma$ such that $H\curvearrowright Y$ is minimal and $H$ is $Y$-plump. Then, every intermediate $C^*$-algebra $\mathcal{A}$ with $C(Y)\rtimes_r H\subseteq \mathcal{A}\subseteq C(X)\rtimes_r\Gamma$, is simple.
\end{maintheorem}
Under additional dynamical assumptions, we also obtain rigidity results for intermediate $C^*$-algebras.
Specifically, we show that in the presence of $X$-plumpness, every intermediate $C^*$-algebra between $C(X)\rtimes_r H$ and $C(X)\rtimes_r \Gamma$ is itself a crossed product $C^*$-algebra. This extends the main rigidity theorem of \cite{amrutam2021intermediate} from the group case to the relative crossed product setting, and yields new classes of examples beyond those previously obtained by Suzuki~\cite{suzuki2020complete} and by the first-named author in~\cite{amrutam2021intermediate, amrutam2023intermediate,amrutam2024crossed}. Below, $\mathbb{E}$ denotes the canonical conditional expectation associated with the reduced crossed product.  
\begin{maintheorem}\label{thmD}
Let $X$ be a $\Gamma$-space, and $\mathcal{B}$ be a unital $\Gamma$-$C^*$-algebra. Consider $C(X)\otimes_{\text{min}}\mathcal{B}$, with the diagonal action.
Let $N = \ker(\Gamma \curvearrowright \mathcal{B})$. Assume that $N$ is $X$-Plump and $N\curvearrowright X$ is minimal. Then for every intermediate C*-algebra $\mathcal{A}$ with
\[
C(X) \rtimes_r \Gamma \subseteq \mathcal{A} \subseteq (C(X) \otimes_{\min} \mathcal{B}) \rtimes_r \Gamma
\]
we have that $E(\mathcal{A}) \subseteq \mathcal{A}$. If, in addition, $\Gamma$ has the approximation property (AP), $\mathcal{A}$ is a crossed product of the form $(C(X) \otimes_{\min} \tilde{\mathcal{B}}) \rtimes_r \Gamma$. Here, $\tilde{\mathcal{B}}$ is a $\Gamma$-invariant $C^*$-subalgebra of $\mathcal{B}$. 
\end{maintheorem}
It is worth emphasizing that such rigidity phenomena are inherently relative and cannot be expected to hold when $X$ is a point.
Indeed, plumpness of a subgroup forces $C^*$-simplicity of the ambient group~\cite{amrutam2021intermediate}. While it is known by the work of Amrutam–Glasner–Glasner~\cite{amrutam2024non} that for non-amenable groups there exist intermediate $C^*$-algebras of the reduced group $C^*$-algebra which are not crossed products.
The presence of a nontrivial base space $X$ thus plays a crucial role in restoring crossed product rigidity.
\stoptocentries
\subsection*{Organization of the Paper}The paper is organized as follows. In Section~\ref{sec:pre}, we review the necessary background on relative
boundaries, strongly proximal extensions, and establish basic stability results for finite-index subgroups. In Section~\ref{sec:uni}, we construct the universal generalized boundary
$RB_X(\Gamma, H)$ using an inverse limit argument, establishing
Theorem~\ref{thmA}.
Section~\ref{sec:rel} focuses on the dynamical extension problem; we prove
the Relative Boundary Extension Theorem (Theorem~\ref{thmB}) and analyze the
freeness of the extended action via the centralizer of the subgroup (see Theorem~\ref{thm:extensionaction}). In this section, we also prove
that for commensurated subgroups, the universal boundary coincides with the
extended generalized Furstenberg boundary (Theorem~\ref{thm:same}). Finally, Section~\ref{sec:application} is devoted to applications.  It is here that we introduce the notion of
$X$-plumpness and prove the simplicity and rigidity results for intermediate
crossed products (Theorems~\ref{thmC} and \ref{thmD}), concluding with
several concrete classes of examples.
\subsection*{Acknowledgments} C. Liu was supported by   the
		Postdoctoral Fellowship Program and China Postdoctoral
		Science Foundation under Grant Number BX20250067, and the China Postdoctoral Science
		Foundation under Grant Number 2025M773074. The authors are grateful to Yongle Jiang for taking the time to read a near-complete draft of this paper and for his several comments and corrections. The first-named author thanks Mehrdad Kalantar for many helpful discussions. We also thank Kang Li, Eduardo Scarparo and Zahra Naghavi for many useful remarks. 
\starttocentries{}

\section{Preliminaries}\label{sec:pre}
\stoptocentries
\subsection{Dynamics}
In this subsection, we establish the dynamical framework that underpins our results. We review the theory of $\Gamma$-spaces and factor maps, focusing on the hierarchy of extensions—proximal, strongly proximal, and minimal—that generalizes the classical study of boundary actions (initially introduced in \cite{SGlasner1975}). A crucial tool in this analysis is the relative version of Glasner’s rigidity lemma, which asserts that the generalized Furstenberg boundary is the unique universal object in its category. We conclude by examining specific constructions that highlight the existence of generalized boundaries beyond trivial product structures.

Let $X$ be a {\em $\Gamma$-space}, i.e.\ $X$ is a compact Hausdorff space under a continuous action $\Gamma$. We say that the action $\Gamma\curvearrowright X$ is \emph{topologically free} if for every $g\in\Gamma\setminus{e}$ the fixed-point set $\Fix_X(g)$ has empty interior in $X$. A closed invariant subset $B\subset X$ is called
\emph{$\Gamma$-minimal} if for every $x\in B$ one has
$\overline{\Gamma x}=B$. When the acting group is clear from context,
we say that $B$ is a \emph{minimal subset} and omit $\Gamma$. Let $\Prob(X)$ denote the space of Borel probability measures on $X$.
For $\mu\in\Prob(X)$, we write $\supp(\mu)$ for the (topological) support of $\mu$.

A map  $\pi_Y:Y\to X$ between two $\Gamma$-spaces is called a {\em factor map} if $\pi_Y$ is a continuous surjection and satisfies $g\circ\pi_Y=\pi_Y\circ g.$

\begin{definition}
   Let $X$ and $Y$ be two $\Gamma$-spaces, and let  $\pi_Y:Y\to X$ be a factor map. Then $\Ypi$ is called a {\em proximal extension} if for any $y,y'\in \pi_Y^{-1}(x)$ for some $x\in X$, there exists a sequence $\{g_i\}_{i\in\mb N}$ such that
   \[\lim_{i\to\infty}g_iy=\lim_{i\to\infty}g_iy'.\]
\end{definition}

\begin{definition}
   Let $X$ and $Y$ be two $\Gamma$-spaces, and let  $\pi_Y:Y\to X$ be a factor map. Then 
   \begin{enumerate}
       \item $\Ypi$ is called a \emph{minimal extension} if $Y$ is minimal.
       \item $\Ypi$ is called a \emph{strongly proximal extension} if for any $\nu\in\Prob(Y)$ with $\supp(\nu)\subset \pi_Y^{-1}(x)$, for some $x\in X$, there exists $y\in Y$ such that $\delta_y\in\overline{\Gamma\nu}$.
   \end{enumerate}
\end{definition}
To establish a hierarchy among extensions, we first verify that the notion of strong proximality indeed strengthens standard proximality.
\begin{proposition}\label{prop:strongly proximal extension implies proximal extension}
       Let $X$ and $Y$ be two $\Gamma$-spaces, and let  $\pi_Y:Y\to X$ be a factor map. If $\Ypi$ is a strongly proximal extension, then it is a proximal extension.
\end{proposition}
\begin{proof}
    Given  any $y,y'\in \pi_Y^{-1}(x)$ for some $x\in X$, one has 
    \[\supp\left(\frac{1}{2}(\delta_y+\delta_{y'})\right)\subset\pi_Y^{-1}(x).\]
   Since $\Ypi$ is strongly proximal, there exist a net $\{t_i\}_{i\in I}\subset \Gamma$ and $y^*\in Y$ such that  
   \[\lim_{i\in I}\frac{1}{2}\left(\delta_{t_iy}+\delta_{t_iy'}\right)=\lim_{i\in I}t_i\frac{1}{2}\left(\delta_y+\delta_{y'}\right)=\delta_{y^*}.\]
Therefore, one has  
\[\lim_{i\in I}t_iy=\lim_{i\in I}t_iy'=y^*.\]
As $(y,y')\in Y\times Y$ is arbitrary, $\Ypi$ is proximal.
\end{proof}
Proximal extensions impose strong structural constraints on the extension space. Specifically, over a minimal base, the contracting nature of proximality forces the uniqueness of the minimal subset. Although the proof is straightforward, we provide it for completeness.
\begin{lemma}\label{lem:unique of minimal set}
    Let $Y$ be a $\Gamma$-space, and let $X$ a minimal $\Gamma$-space. Suppose that $\pi_Y:Y\to X$ is a proximal extension. Then 
    $Y$ has a unique minimal subset.
\end{lemma}
\begin{proof}
    Suppose, for a contradiction, that there exist two different minimal subsets $Y_1$ and $Y_2$. Since $\pi_Y$ is surjective and $\pi_Y(Y_i)$, $i=1,2$ are closed and invariant, it follows that 
    \[\pi_Y(Y_1)=X=\pi_Y(Y_2).\]
    Choose $y_i\in Y$, $i=1,2$ such that $\pi(y_1)=\pi(y_2)$. Since $\pi_Y$ is proximal, it follows that 
    \[\lim_{i\to\infty}g_iy=\lim_{i\to\infty}g_iy'.\]
    This implies that $Y_1=Y_2$, a contradiction. 
\end{proof}
A central theme of this work is the relationship between the dynamics of a group and its subgroups. We first observe that proximality passes to finite-index subgroups.
\begin{lemma}\label{lem:proximal preserved by normal subgroups}
    	Let $X$ be a minimal $H$-space and let $K\leq H$ be a 
		subgroup with $[K:H]<\infty$. Let $\pi_Y:Y\to X$ be a factor map.
        If $\pY$ is $H$-proximal, then $\pY$ is $K$-proximal.
\end{lemma}
\begin{proof}
    Since $\pY$ is $H$-proximal, one has  for any $(y,y')\in Y\times Y$, there exists $\{h_i\}_{i\in\mb N}$ such that 
\[\lim_{i\to\infty}h_iy=\lim_{i\to\infty}h_iy'.\]
Since $[K:H]<\infty$, there exist $h\in H$ and a subsequence $\{k_i\}_{i\in\mb N}$ such that 
\[\lim_{i\to\infty}hk_iy=\lim_{i\to\infty}hk_iy'.\]
By continuity, passing to a subsequence, it follows that 
\[\lim_{i\to\infty}k_iy=\lim_{i\to\infty}k_iy',\]
which implies that $\pY$ is $K$-proximal, as $(y,y')\in Y\times Y$ is arbitrary.
\end{proof}
Similarly, minimality behaves well under finite extensions, provided we restrict our attention to normal subgroups and proximal spaces. This allows us to descend dynamical properties from $H$ to its finite-index subgroups.
\begin{lemma}\label{lem:minimal preserved by normal subgroups}
    	Let $X$ be a minimal $H$-space and let $K\leq H$ be a normal
		subgroup with $[H:K]<\infty$. Let $\pi_Y:Y\to X$ be a factor map.
        If $\pY$ is $H$-proximal and $Y$ is $H$-minimal, then $Y$ is $K$-minimal.
\end{lemma}
\begin{proof}
Since $\pY$ is $H$-proximal, it follows from Lemma \ref{lem:proximal preserved by normal subgroups} that $\pY$ is $K$-proximal. Combining this with Lemma \ref{lem:unique of minimal set}, one has $Y$ has a unique $K$-minimal subset, denoted by $M_Y$. 

We now prove that $Y=M_Y$, and hence $Y$ is $K$-minimal. To that end, we consider the set $hM_Y$ for all $h\in H$. Since $K$ is a normal subgroup of $H$, for any $k\in K$, there exists $k'\in K$ such that 
$kh=hk'$. Since $M_Y$ is $K$-invariant, it follows that 
\[k(hM_Y)=hk'M_Y=hM_Y.\]
Therefore, $hM_Y\subset M_Y$ for all $h\in H$. Thus, $M_Y$ is an $H$-invariant closed subset of $Y$, which implies that $M_Y=Y$.
\end{proof}
We recall the following definitions from \cite{SGlasner1975}.
\begin{definition}\label{def:GX-boundary}
Let $X$ be a $\Gamma$-space. A \emph{$(\Gamma,X)$-boundary} is a pair $(Y,\pi_Y)$ such that $\Ypi$ is a minimal strongly proximal extension of $X$.
\end{definition}

\begin{definition}[Generalized Furstenberg boundary]\label{def:gen-F-boundary}
 A $(\Gamma,X)$-boundary
$
  (\partial_F(\Gamma,X),\pi_\Gamma)
$
is called the \emph{generalized Furstenberg boundary} of $(\Gamma,X)$ if it
has the following universal property: for every $(\Gamma,X)$-boundary
$(Y,\pi_Y)$ there exists a unique continuous $\Gamma$-equivariant surjection
\[
  \Phi:\partial_F(\Gamma,X)\to Y
\]
such that the diagram
\[
\begin{tikzcd}
\partial_F(\Gamma,X) \arrow[dr, "\pi_\Gamma"'] \arrow[rr, "\Phi"] & & 
Y \arrow[dl, "\pi_Y"] \\
& X &
\end{tikzcd}
\]
commutes, i.e.\ $\pi_Y\circ \Phi = \pi_\Gamma$.
\end{definition}

\begin{remark}
The generalized Furstenberg boundary exists and is unique up to $\Gamma$-isomorphism in the sense that if $(Y,\pi_Y)$ is another $(\Gamma,X)$-boundary with the same universal
property, then there exists a unique $G$-isomorphism
$\Psi: \partial_F(\Gamma,X)\to Y$ satisfying $\pi_Y\circ\Psi=\pi_G$ (see \cite[p. 7]{naghavi2020furstenberg}).

In the special case where $X$ consists of a single point (with the
trivial $\Gamma$-action), $(\partial_F(\Gamma,X),\pi_\Gamma)$ reduces to the classical
Furstenberg boundary $\partial_F \Gamma$.
\end{remark}
A defining feature of the classical Furstenberg boundary is its rigidity: the only equivariant self-map is the identity. The following lemma generalizes this fundamental rigidity principle to the relative setting of extensions. It is a relative version of \cite[Chapter II, Lemma 4.1]{SGlasnerbook}. 
\begin{lemma}\label{lem:relative-Glasner}
Let $X$ be a $\Gamma$-space and let $(Y,\pi_Y)$, $(Z,\pi_Z)$ be two
$(\Gamma,X)$-boundaries. Suppose $f,g:Y\to Z$ are continuous
$\Gamma$-equivariant maps such that
\begin{equation}\label{eq:self-homo}
    \pi_Z\circ f = \pi_Y = \pi_Z\circ g.
\end{equation}
Then $f=g$, i.e.\ the following diagram commute:
\[\begin{tikzcd}
Y \arrow[dr, "\pi_Y"'] \arrow[rr, shift left=0.6ex, "f"] 
                      \arrow[rr, shift right=0.6ex, "g"'] & &
Z \arrow[dl, "\pi_Z"] \\
& X &
\end{tikzcd}\]

In particular, any $\Gamma$-equivariant continuous self-map $u$ of
$\partial_F(\Gamma,X)$ with $\pi_\Gamma\circ u=\pi_\Gamma$ is the identity,
i.e.,\ the following diagram commutes.
\[
\begin{tikzcd}
\partial_F(\Gamma,X) \arrow[rr, "u"] \arrow[dr, "\pi_\Gamma"'] & &
\partial_F(\Gamma,X) \arrow[dl, "\pi_\Gamma"] \\
& X &
\end{tikzcd}
\]
\end{lemma}
\begin{proof}
Set
\[
   E := \{y\in Y : f(y)=g(y)\}.
\]
We show that $E$ is non-empty, closed, and $\Gamma$-invariant, and hence $E=Y$, as $Y$ is minimal.

Since  $f$ and $g$ are continuous, $E$ is closed. Meanwhile, as $f,g$ are $\Gamma$-equivariant, it follows that $E$ is $\Gamma$-equivariant.
Finally, to see that $E$ is non-empty, fix $y_0\in Y$ and consider the measure
\[
   \mu := \frac12\bigl(\delta_{f(y_0)}+\delta_{g(y_0)}\bigr)
   \in\Prob(Z).
\]
Note that 
\[\pi_Z(f(y_0))=\pi_Y(y_0)=\pi_Z(g(y_0)),\]
and hence
\[\supp(\mu)\subset \pi_Z^{-1}(\pi_Y(y_0)).\]
Since $(Z,\pi_Z)$ is a strongly proximal extension, there exist a net $(t_i)_{i\in I}$
in $\Gamma$ and a point $z\in Z$ such that $\lim_it_i\mu=\delta_z$ in
$\Prob(Z)$. 

On the other hand,
\[
   t_i\mu
   = \frac12\bigl(\delta_{t_i f(y_0)}+\delta_{t_i g(y_0)}\bigr)
   = \frac12\bigl(\delta_{f(t_i y_0)}+\delta_{g(t_i y_0)}\bigr),
\]
where we have used $\Gamma$-equivariance of $f$ and $g$ in the last equality.
Thus, one has 
\[
\lim_ i f(t_i y_0)=z,\text{ and }\lim_i g(t_i y_0)= z.
\]
By compactness of $Y$ there exists a subnet $(t_{i_j}y_0)_{j}$ converging
to some $y_*\in Y$. By continuity,
\[
   f(y_*) = z = g(y_*),
\]
which implies that $y_*\in E$. Therefore, $E\neq \emptyset$.

By the discussion above, we must have $E=Y$. Thus $f(y)=g(y)$ for all
$y\in Y$, i.e.\ $f=g$.

For the particular case, let $Y=Z=\partial_F(\Gamma,X)$,
$\pi_Y=\pi_Z=\pi_\Gamma$, and let $u:\partial_F(\Gamma,X)\to\partial_F(\Gamma,X)$
be a $\Gamma$-equivariant self-homeomorphism with $\pi_\Gamma\circ u=\pi_\Gamma$. By the above argument with $f=u$ and $g=\mathrm{Id}_{\partial_F(\Gamma,X)}$,
one has $u=\mathrm{Id}_{\partial_F(\Gamma,X)}$.
\end{proof}
Despite the theoretical importance of generalized Furstenberg boundaries, to the best of our knowledge, there are relatively few explicit examples. The following constructions address this by providing a natural class of examples arising from diagonal actions.
\begin{example}[Canonical boundary extension]\label{Ex:strongproximal}
Let $\Gamma$ be a non-amenable group and let $X$ be a minimal $\Gamma$-space.
Let $B$ be a non-trivial $\Gamma$-boundary.
We now construct a $(\Gamma,X)$-boundary in a canonical way. Consider the product space
\[
   Y := B\times X
\]
equipped with the diagonal $\Gamma$-action, and let
\[
   \pi:Y\longrightarrow X,\qquad \pi(b,x)=x
\]
be the projection onto the second coordinate.
We claim that $\pi$ is a strongly proximal extension. Let $\nu\in \Prob(Y)$ be such that
\[
   \supp(\nu)\subset \pi^{-1}(x)
   \quad\text{for some }x\in X.
\]
Since $\pi^{-1}(x)=B\times\{x\}$, the measure $\nu$ is supported on
$B\times\{x\}$ and hence can be written as
\[
   \nu = \eta\otimes \delta_x
\]
for some $\eta\in \Prob(B)$.
Because $\Gamma\curvearrowright B$ is a boundary action, there exists a net
$(s_j)_j\subset \Gamma$ and a point $b\in B$ such that
\[
   s_j\eta \longrightarrow \delta_b
   \quad\text{in }\Prob(B).
\]
Passing to a subnet if necessary, we may also assume that $s_jx\to x'$ for some
$x'\in X$.
Then
\[
   s_j\nu
   = (s_j\eta)\otimes \delta_{s_jx}
   \longrightarrow \delta_b\otimes \delta_{x'}
\]
which is a Dirac measure supported on the fibre $\pi^{-1}(x')$.
This shows that $\pi$ is a strongly proximal extension.

By Proposition~\ref{prop:strongly proximal extension implies proximal extension},
$\pi$ is proximal.
Since the base system $X$ is minimal, Lemma~\ref{lem:unique of minimal set}
implies that $Y$ contains a unique minimal $\Gamma$-invariant subset; denote it
by $\hat X$.
The restriction
\[
   \pi|_{\hat X} : \hat X \longrightarrow X
\]
is then a minimal, strongly proximal extension.
In particular, $(\hat X,\pi|_{\hat X})$ is a $(\Gamma,X)$-boundary.
\end{example}
\begin{remark}
In general, the $(\Gamma,X)$-boundary $\hat X$
constructed in Example~\ref{Ex:strongproximal} does not coincide with
the generalized Furstenberg boundary of $(\Gamma,X)$. Let $\Gamma$ be a non-amenable group. By \cite[Theorem~III.6.1]{SGlasnerbook}, every minimal strongly proximal
$\Gamma$-flow is disjoint from every minimal $\Gamma$-flow that admits
an invariant probability measure. In particular, let $X$ be an infinite minimal
$\Gamma$-space carrying a $\Gamma$-invariant Borel probability measure.
Then $X$ and $\partial_F\Gamma$ are disjoint, and hence the diagonal
action $\Gamma \curvearrowright X\times\partial_F\Gamma$ is minimal.
In this situation, Example~\ref{Ex:strongproximal} applied to
$B=\partial_F\Gamma$ shows that
\[
   \hat X = \partial_F\Gamma\times X
\]
is a generalized $(\Gamma,X)$-boundary.

Suppose towards a contradiction that $\hat X$ is (isomorphic to) the
generalized Furstenberg boundary $\partial_F(\Gamma,X)$. By \cite[Proposition 3.3]{kawabe2017uniformlyrecurrentsubgroupsideal}, $\partial_F(\Gamma,X)$ is  Stonean. Hence $\hat X$ would have to be Stonean. Since $\Gamma$ is non-amenable, it follows from \cite[Corollary 3.13]{extremallytopologicalextremallyspace} that $X$ is finite. This is a contradiction. So, $\hat{X}$ can not be isomorphic to the generalized Furstenberg boundary in general.
\end{remark}
Of course, the above construction does not work for amenable groups. However, in the case when the group has a non-trivial amenable radical, it is possible to construct such boundaries over any given space.
\begin{example}\label{Ex:productgroup}
Let $\Gamma$ be a non-amenable group with the non-trivial amenable radical
$N$.
Let $X$ be a compact $\Gamma$-space such that the restricted action
$N\curvearrowright X$ is minimal.
Denote by $\partial_F\Gamma$ the Furstenberg boundary of $\Gamma$.

By \cite[Corollary~8]{furman2003minimal}, the amenable radical $N$ is
contained in the kernel of the action on $\partial_F\Gamma$, i.e.
$N\subset \ker(\Gamma\curvearrowright\partial_F\Gamma)$.
Since every $\Gamma$-boundary $B$ is a $\Gamma$-factor of
$\partial_F\Gamma$, it follows that $N$ acts trivially on $B$ as well.

We claim that for any $\Gamma$-boundary $B$, the pair
\[
   (X\times B,\ \pi_X),
   \qquad \pi_X(x,b)=x,
\]
is a $(\Gamma,X)$-boundary.

First, as in Example~\ref{Ex:strongproximal}, consider the product
$Y:=X\times B$ with the diagonal $\Gamma$-action and the projection
$\pi_X:Y\to X$.
The same argument as in Example~\ref{Ex:strongproximal} shows that
$\pi_X$ is a strongly proximal extension:
whenever $\nu\in\Prob(Y)$ is supported on a fibre $\pi_X^{-1}(x)$, one
can find a net $(\gamma_i)_i\subset\Gamma$ such that
$\gamma_i\nu\to \delta_{(x',b')}$ for some $(x',b')\in Y$, and hence
$\pi_X$ is strongly proximal over $X$.

It remains to show that the action $\Gamma\curvearrowright Y$ is minimal.
Fix $(x,b)\in X\times B$ and let $(x',b')\in X\times B$ be arbitrary.
We will show that $(x',b')$ lies in the closure of the $\Gamma$-orbit of
$(x,b)$.

Since $\Gamma\curvearrowright B$ is minimal and $B$ is a
$\Gamma$-boundary, there exists a net $(\gamma_i)_i\subset\Gamma$ such
that
\[
   \gamma_i b \longrightarrow b'.
\]
Passing to a subnet if necessary, we may assume that
\[
   \gamma_i x \longrightarrow x_0
\]
for some $x_0\in X$.
By minimality of the $N$-action on $X$, there exists a net
$(n_j)_j\subset N$ such that
\[
   n_j x_0 \longrightarrow x'.
\]

Consider now the net $(n_j\gamma_i)_{(j,i)}$ in $\Gamma$ indexed by the
directed set of pairs.
Then
\[
   (n_j\gamma_i)x = n_j(\gamma_i x)\longrightarrow x',
\]
and, since $N$ acts trivially on $B$,
\[
   (n_j\gamma_i)b = n_j(\gamma_i b) = \gamma_i b \longrightarrow b'.
\]
Thus
\[
   (n_j\gamma_i)(x,b)\longrightarrow (x',b'),
\]
showing that the $\Gamma$-orbit of $(x,b)$ is dense in $X\times B$.
Since $(x,b)$ was arbitrary, the action $\Gamma\curvearrowright X\times B$
is minimal.

Summarizing, $\pi_X:X\times B\to X$ is a minimal strongly proximal
extension, so $(X\times B,\pi_X)$ is a $(\Gamma,X)$-boundary.
\end{example}

To patch together the local actions defined by commensuration, we need to transport boundary structures via isomorphisms. The next lemma allows us to canonically identify boundaries of isomorphic groups acting on equivalent spaces.
\begin{lemma}
    \label{lem:rel-iso-boundary}
Let $\Gamma_1,\Gamma_2$ be countable discrete groups and let
		$\psi:\Gamma_1\to\Gamma_2$ be a group isomorphism.
		Let $X_1$ be a $\Gamma_1$-space and $X_2$ a $\Gamma_2$-space.
		Assume there exists a homeomorphism
		$\phi:X_1\to X_2$
		such that
		\begin{equation}\label{eq:rel-iso-base}
			\phi(gx)=\psi(g)\,\phi(x)
			\qquad\text{for all }g\in\Gamma_1,\ x\in X_1 .
		\end{equation}
		Then there exists a unique homeomorphism
		\[
		\widetilde\psi:\partial_F(\Gamma_1,X_1)\longrightarrow
		\partial_F(\Gamma_2,X_2)
		\]
		such that
		\begin{enumerate}
			\item $\widetilde\psi$ is $\Gamma_1$-equivariant in the sense that
			\[
			\widetilde\psi(g y)=\psi(g)\,\widetilde\psi(y),
			\qquad \forall g\in\Gamma_1,\ y\in\partial_F(\Gamma_1,X_1),
			\]
			where $\Gamma_1$ acts on $\partial_F(\Gamma_2,X_2)$ via $\psi$;
			\item the factor maps satisfy
			\[
			\phi\circ \pi_{\Gamma_1}
			\;=\; \pi_{\Gamma_2}\circ\widetilde\psi,
			\]
			i.e.\ the following diagram commutes:
			\[
			\begin{tikzcd}
				\partial_F(\Gamma_1,X_1) \arrow{r}{\widetilde\psi}
				\arrow{d}[swap]{\pi_{\Gamma_1}} &
				\partial_F(\Gamma_2,X_2) \arrow{d}{\pi_{\Gamma_2}}\\
				X_1 \arrow{r}{\phi} & X_2.
			\end{tikzcd}
			\]
		\end{enumerate}
\end{lemma}
\begin{proof}
Consider the action $\Gamma_1\curvearrowright \partial_F(\Gamma_2,X_2)$ defined by
\[
   g\cdot y := \psi(g)y,\qquad g\in \Gamma_1,\ y\in \partial_F(\Gamma_2,X_2).
\]
Thus $\partial_F(\Gamma_2,X_2)$ becomes a $\Gamma_1$-space.  
Since $\pi_{\Gamma_2}:\partial_F(\Gamma_2,X_2)\to X_2$ is $\Gamma_2$-equivariant, we define
\[
   \pi := \phi^{-1}\circ \pi_{\Gamma_2}:\partial_F(\Gamma_2,X_2)\to X_1.
\]
Then for any $g\in\Gamma_1$ and $y\in \partial_F(\Gamma_2,X_2)$ we have
\[
\begin{aligned}
   \pi(g\cdot y)
   &= \phi^{-1}\bigl(\pi_{\Gamma_2}(\psi(g)\,y)\bigr)
    = \phi^{-1}\bigl(\psi(g)\,\pi_{\Gamma_2}(y)\bigr) \\
   &= g \phi^{-1}\bigl(\pi_{\Gamma_2}(y)\bigr)
    = g \pi(y),
\end{aligned}
\]
showing that $\pi$ is $\Gamma_1$-equivariant.  
Since $\psi$ is a group isomorphism and the $\Gamma_2$-action on $\partial_F(\Gamma_2,X_2)$ is minimal and strongly proximal, the induced $\Gamma_1$-action is also minimal and strongly proximal. Hence $\pi$ is precisely the factor map making $\partial_F(\Gamma_2,X_2)$ into a $(\Gamma_1,X_1)$-boundary.

The
universal property yields a unique continuous $\Gamma_1$-equivariant
surjection
\[
   \widetilde\psi: \partial_F(\Gamma_1,X_1)\to
   \partial_F(\Gamma_2,X_2)
\]
such that
\[
   \pi\circ \widetilde\psi = \pi_{\Gamma_1} \implies  \pi_{\Gamma_2}\circ \widetilde\psi =\phi\circ \pi_{\Gamma_1}.
\]
Thus, this map satisfies (i) and (ii) by construction.

Therefore, it suffices to prove that $\widetilde\psi$ is a homeomorphism.
Applying the same construction with the inverse isomorphism
\[
   \psi^{-1}: \Gamma_2\to\Gamma_1,
\]
we regard $\partial_F(\Gamma_1,X_1)$ as a $(\Gamma_2,X_2)$-boundary
and obtain a continuous $\Gamma_2$-equivariant surjection
\[
   \widehat\psi: \partial_F(\Gamma_2,X_2)\to
   \partial_F(\Gamma_1,X_1)
\]
such that
\[
   \pi_{\Gamma_1}\circ \widehat\psi =\phi^{-1} \circ \pi_{\Gamma_2}.
\]

Consider the compositions
\[
   \alpha := \widehat\psi\circ\widetilde\psi
   : \partial_F(\Gamma_1,X_1)\to\partial_F(\Gamma_1,X_1),
\]
\[
   \beta := \widetilde\psi\circ\widehat\psi
   : \partial_F(\Gamma_2,X_2)\to\partial_F(\Gamma_2,X_2).
\]
Since $\widetilde\psi$ satisfies (ii), it follows that $\alpha$ is $\Gamma_1$-equivariant and
\[
   \pi_{\Gamma_1}\circ\alpha
   = \pi_{\Gamma_1}\circ\widehat\psi\circ\widetilde\psi
   = \phi^{-1}\circ \pi_{\Gamma_2}\circ\widetilde\psi
   = \pi_{\Gamma_1}.
\]
Similarly, $\beta$ is $\Gamma_2$-equivariant and
\[
   \pi_{\Gamma_2}\circ\alpha
   = \pi_{\Gamma_2}.
\]

By Lemma~\ref{lem:relative-Glasner} applied
to $(\Gamma_1,X_1)$ and $(\Gamma_2,X_2)$ respectively, both $\alpha$ and
$\beta$ must be the identity maps:
\[
   \alpha = \mathrm{id}_{\partial_F(\Gamma_1,X_1)},\qquad
   \beta = \mathrm{id}_{\partial_F(\Gamma_2,X_2)}.
\]
Thus $\widehat\psi$ is the inverse of $\widetilde\psi$, which is
therefore a homeomorphism. 
\end{proof}

The key step in promoting a relative boundary action from a commensurated subgroup to the whole group is to understand how the relative boundary behaves upon restriction. The following lemma establishes that the generalized boundary is stable—and effectively invariant—under passage to finite-index subgroups.
\begin{lemma}\label{lem:key lemma}
    	Let $X$ be a $H$-space and let $K\leq H$ be a
		subgroup of finite index, $[H:K]<\infty$. Let $\pi_Y:Y\to X$ be a factor map.
        Then we have the following:
        \begin{enumerate}
            \item If $(Y,\pi_Y)$ is  a strongly proximal extension as an $H$-space, then it is  a strongly proximal extension as a $K$-space.
            \item If $(Y,\pi_Y)$ is  an $(H,X)$-boundary, then it is a $(K,X)$-boundary.
            \end{enumerate}
\end{lemma}
\begin{proof}
     Since  $[H:K]<\infty$, there exist $h_1,\ldots,h_n$ for some $n\in\N$ such that 
  \begin{equation}\label{eq:finite co}
      H=\bigsqcup_{j=1}^nh_j K.
  \end{equation}  
    
\noindent\textbf{(i):}  For  any $\nu\in \Prob(Y)$ with $\supp(\nu)\subset\pi^{-1}_Y(x)$ for some $x\in X$, as $(H,Y)$ is strongly proximal, it follows that there exist $y\in Y$ and a net $\{t_i\}_{i\in I}\subset I$ such that $\lim_it_i \nu=\delta_y$.  By \eqref{eq:finite co}, passing to a subsequence, we may assume, without loss of generality, that there exists $j\in\{1,2,\ldots,n\}$ and  a net $\{k_i\}_{i\in I}\subset K$ such that
\[t_i=h_j k_i,\text{ for all }i\in I.\]
Then
\[\delta_y=\lim_it_i \nu=\lim_ih_jk_i \nu,\]
which implies that 
\[\lim_ik_i \nu=(h_j^{-1})_*\delta_y=\delta_{h_j^{-1}y}.\]
Since $h_j^{-1}y\in Y$, and $\nu\in\Prob(Y)$ is arbitrary, it follows that $(Y,\pi_Y)$ is  a strongly proximal extension as a $K$-space.

\medskip

\noindent\textbf{(ii):} Since $(Y,\pi_Y)$ is  an $(H,X)$-boundary, one has $\Ypi$ is a minimal strongly proximal extension as an $H$-space. By (i), 
we only need to prove that $Y$ is  minimal as a $K$-space. 

Since $[H:K]<\infty$, it follows from \cite[p. 26]{HewitRoss1979book} that there exists a normal subgroup $N$ of $H$ which is contained in $K$ with $[H:N]<\infty$ such that $N\subset K$. We now know that $\Ypi$ is a strongly proximal extension, it follows from Proposition \ref{prop:strongly proximal extension implies proximal extension} that it is a proximal extension. By Lemma \ref{lem:proximal preserved by normal subgroups} and Lemma \ref{lem:minimal preserved by normal subgroups}, one has 
$Y$ is  minimal as a $N$-space. Since $N\subset K$, it follows that $Y$ is minimal as a $K$-space. 

Combining these two facts, we have $\Ypi$  is a $(K,X)$-boundary.
\end{proof}
While the above lemma ensures that boundary properties descend to finite-index subgroups, we must also understand how they interact with the ambient group structure. In the case of a normal subgroup, the equivariance of the boundary map forces a strong rigidity in the sense that the action of the normal subgroup uniquely determines an action of the entire group. We follow the proof in \cite[Proposition~4.3]{SGlasnerbook} and adapt it to our setup.
\begin{proposition}
\label{prop:normalextension}
Let $N\triangleleft\Gamma$ be a normal subgroup and  let $X$ be a $\Gamma$-space
		which is minimal as an $N$-space. Then, $N\curvearrowright\partial_F(N,X)$ can be extended to $\Gamma\curvearrowright\partial_F(N,X)$. 
\end{proposition}
\begin{proof} For each $h\in\Gamma$, define the conjugation automorphism
		\[
		\psi_h:N\to N,\qquad \psi_h(n):=hnh^{-1}.
		\]
		Since $N\triangleleft\Gamma$, $\psi_h$ is a group automorphism of $N$.
		Also define a homeomorphism $\phi_h:X\to X$ by $\phi_h(x):=h x$.
		Then for every $n\in N$ and $x\in X$ we have
		\[
		\phi_h(n x)=h(nx)=(hnh^{-1})(h x)=\psi_h(n)\phi_h(x).
		\]
		Thus $\phi_h$ satisfies \eqref{eq:rel-iso-base} in Lemma~\ref{lem:rel-iso-boundary}
		with $\Gamma_1=\Gamma_2=N$, $X_1=X_2=X$, and $\psi=\psi_h$. Therefore, applying Lemma~\ref{lem:rel-iso-boundary}, we obtain a unique homeomorphism
		\[
		T_h:\partial_F(N,X)\longrightarrow \partial_F(N,X)
		\]
		such that
		\begin{align}
			&T_h(n y)=\psi_h(n) T_h(y),
			 &\forall\,n\in N,\ y\in \partial_F(N,X), \label{eq:Th-conj}\\
		& \pi_N(T_h(y))=h\pi_N(y)=\varphi_h(\pi_N(y))\ \ &\forall y\in \partial_F(N,X). \label{eq:Th-cover}
		\end{align}
In particular, we have the following commutative diagram.
\[
			\begin{tikzcd}
				\partial_F(N,X) \arrow{r}{T_h}
				\arrow{d}[swap]{\pi_{N}} &
				\partial_F(N,X) \arrow{d}{\pi_{N}}\\
				X \arrow{r}{\phi_h} & X.
			\end{tikzcd}
			\]
Now define an action of $\Gamma$ on $\partial_F(N,X)$ by
		\[
		hy:=T_h(y),\qquad h\in\Gamma,\ y\in \partial_F(N,X).
		\]
We have the following two claims.	\begin{claim}\label{cl1}
	$T_{h_1h_2}=T_{h_1}\circ T_{h_2}$ for all $h_1,h_2\in\Gamma$.
\end{claim}
\begin{proof}[Proof of Claim \ref{cl1}]
Using $\phi_{h_1h_2}=\phi_{h_1}\circ\phi_{h_2}$, we see that
\[
\pi_N\circ(T_{h_1}\circ T_{h_2})
\stackrel{eqn~\eqref{eq:Th-cover}}{=}(\phi_{h_1}\circ\pi_N)\circ T_{h_2}
\stackrel{eqn~\eqref{eq:Th-cover}}{=}\phi_{h_1}\circ(\phi_{h_2}\circ\pi_N)
=\phi_{h_1h_2}\circ\pi_N
\stackrel{eqn~\eqref{eq:Th-cover}}{=}h_1h_2\pi_N.
\]
Furthermore, using the identity $\psi_{h_1h_2}=\psi_{h_1}\circ\psi_{h_2}$,
 we have that  for all $n\in N$ and $y\in \partial_F(N,X)$,
\begin{align*}
(T_{h_1}\circ T_{h_2})(n y)
&\stackrel{eqn~\eqref{eq:Th-conj}}{=} T_{h_1}(\psi_{h_2}(n) T_{h_2}(y))
\\&\stackrel{eqn~\eqref{eq:Th-conj}}{=}\psi_{h_1}(\psi_{h_2}(n)) (T_{h_1}\circ T_{h_2})(y)
=\psi_{h_1h_2}(n) (T_{h_1}\circ T_{h_2})(y).
\end{align*}
Hence $T_{h_1}\circ T_{h_2}$ satisfies the same two properties
\eqref{eq:Th-conj}--\eqref{eq:Th-cover} as $T_{h_1h_2}$.
By the uniqueness part of Lemma~\ref{lem:rel-iso-boundary}, we get
$T_{h_1h_2}=T_{h_1}\circ T_{h_2}$. Since $T_e=\text{id}$, we see that \[\text{id}=T_e=T_{hh^{-1}}=T_h\circ T_{h^{-1}},~\forall h\in\Gamma.\]
Consequently, we see that $T_h^{-1}=T_{h^{-1}}$ for all $h\in\Gamma$.
\end{proof}
	
\begin{claim}\label{cl2}
The $\Gamma$-action extends the original $N$-action.
\end{claim} 
\begin{proof}[Proof of Claim \ref{cl2}]
	For $h\in N$, we denote the original action by $L_h$ defined by $L_h(y):=h y$.
Then $\pi_N\circ L_h = h\pi_N$, so $L_h$ satisfies equation~\eqref{eq:Th-cover}.
Also for any $n\in N$ and $y\in \partial_F(N,X)$,
\[
L_h(n y)=h n  y=(h n h^{-1})  (h y)=\psi_h(n)  L_h(y),
\]
so $L_h$ satisfies equation~\eqref{eq:Th-conj}.
By uniqueness in Lemma~\ref{lem:rel-iso-boundary}, $T_h=L_h$, i.e. $T_h(y)=h y$.
\end{proof}

		Combining Claim~\ref{cl1} and Claim~\ref{cl2}, the rule $(h,y)\mapsto T_h(y)$ defines a group action
		of $\Gamma$ on $Y$ by homeomorphisms, and it restricts to the given $N$-action.
		Finally, \eqref{eq:Th-cover} shows that $\pi_N$ is $\Gamma$-equivariant.
	\end{proof}
This extension result enables us to compare boundaries across different subgroups by transforming them to a common normal core. By combining the restriction stability (Lemma~\ref{lem:key lemma}) with the universality of the extended action~Proposition~\ref{prop:normalextension}, we can now show that the generalized Furstenberg boundary is an invariant up to finite index.
\begin{theorem}
\label{thm:equaluptofiniteindex}
Let $K\le H$ be a finite index subgroup, and $X$ a minimal $K$-space.
Both $\partial_F(H,X)$ and $\partial_F(K,X)$ are $(K,X)$-boundaries. Furthermore, there exists  a unique $K$-equivariant homeomorphism
\[
   \varphi: \partial_F(H,X)\to\partial_F(K,X)
\]
such that $\pi_K\circ \varphi = \pi_H$. In particular,
$\partial_F(H,X)$ and $\partial_F(K,X)$ are $K$-isomorphic as
$(K,X)$-boundaries.        
\end{theorem}
\begin{proof}
Let $N$ be a normal subgroup of $H$ contained in $K$ which still has finite index inside $H$ (see e.g. \cite[p. 26]{HewitRoss1979book}).
By the universal property of $\pF(H,X)$, using Proposition~\ref{prop:normalextension}, there exists a  $H$–factor map
\begin{equation}\label{eq:alpha-def}
    \alpha:\pF(H,X)\longrightarrow \pF(N,X)
    \quad\text{with}\quad
    \pi_N\circ\alpha=\pi_H .
\end{equation}
Using Lemma~\ref{lem:key lemma}~(ii), both $\pF(H,X)$ and $\pF(K,X)$ are $(K,X)$–boundaries.  
Therefore, by universality again, there exists a $K$–equivariant factor
map
\begin{equation}\label{eq:psi-def}
   \psi:\pF(K,X)\longrightarrow \pF(H,X)
   \quad\text{such that}\quad
   \pi_H\circ\psi=\pi_K ,
\end{equation}
where $\pi_H:\pF(H,X)\to X$ and $\pi_K:\pF(K,X)\to X$ are the respective
factor maps.  
Similarly, viewing $\pF(N,X)$ and $\pF(K,X)$ as $(N,X)$–boundaries, we
obtain an $N$–equivariant factor map
\begin{equation}\label{eq:eta-def}
   \eta:\pF(N,X)\longrightarrow \pF(K,X)
   \quad\text{with}\quad
   \pi_K\circ\eta=\pi_N .
\end{equation}
Consider now the composition
\[
   \pF(N,X)\xrightarrow{\ \eta\ }\pF(K,X)
   \xrightarrow{\ \psi\ }\pF(H,X)
   \xrightarrow{\ \alpha\ }\pF(N,X).
\]
This is a continuous $N$–equivariant self–map of $\pF(N,X)$, and using
\eqref{eq:alpha-def}–\eqref{eq:eta-def} we compute
\[
   \pi_N\circ\alpha\circ\psi\circ\eta
   = \pi_H\circ\psi\circ\eta
   = \pi_K\circ\eta
   = \pi_N.
\]
By Lemma~\ref{lem:relative-Glasner}, any $N$–equivariant self–map of
$\pF(N,X)$ that fixes $\pi_N$ must be the identity. Hence
\begin{equation}\label{eq:comp-id}
   \alpha\circ\psi\circ\eta
   = \mathrm{id}_{\pF(N,X)},
\end{equation}
which together with the fact that $\al$, $\psi$ and $\eta$ are surjective, implies that $\psi$ is homeomorphism.

Finally, $\psi$ is $K$–equivariant by construction, so
$\phi:=\psi^{-1}:\pF(H,X)\to\pF(K,X)$ is a $K$–equivariant
homeomorphism with $\pi_K\circ\phi=\pi_H$.  
This shows that $\pF(H,X)$ and $\pF(K,X)$ are $K$–isomorphic as
$(K,X)$–boundaries, completing the proof.
\end{proof}

\subsection{Crossed product $C^*$-algebra}
\label{subsec:crossed}
We briefly review the construction of the reduced crossed product, referring the reader to \cite{BO:book} for comprehensive details. Let $\mathcal{A}$ be a unital $\Gamma$-$C^*$-algebra, where the discrete group $\Gamma$ acts on $\mathcal{A}$ by $*$-automorphisms (denoted by $\alpha$).

Fix a faithful $*$-representation $\pi:\mathcal{A} \to \mathbb{B}(\mathcal{H})$. We consider the Hilbert space of square-summable $\mathcal{H}$-valued functions on $\Gamma$:
\[
\ell^2(\Gamma,\mathcal{H}) = \left\{\xi:\Gamma\to \mathcal{H} : \sum_{t\in\Gamma}\|\xi(t)\|_{\mathcal{H}}^2 < \infty \right\}.
\]
We define a unitary representation $\lambda: \Gamma \to \mathbb{B}(\ell^2(\Gamma,\mathcal{H}))$ acting by left translation:
\[
(\lambda_s\xi)(t) = \xi(s^{-1}t), \quad \text{for } s,t \in \Gamma.
\]
Additionally, we define a faithful $*$-representation $\sigma:\mathcal{A} \to \mathbb{B}(\ell^2(\Gamma,\mathcal{H}))$ by
\[
(\sigma(a)\xi)(t) = \pi(\alpha_{t^{-1}}(a))\xi(t), \quad \text{for } a \in \mathcal{A}.
\]
The \emph{reduced crossed product $C^*$-algebra}, denoted by $\mathcal{A}\rtimes_{r}\Gamma$, is defined as the norm-closure of the span of the twisted products inside $\mathbb{B}(\ell^2(\Gamma,\mathcal{H}))$:
\[
\mathcal{A}\rtimes_{r}\Gamma = \overline{\text{Span}\left\{\sigma(a)\lambda_s : a\in\mathcal{A}, s\in\Gamma\right\}}^{\|\cdot\|}.
\]
The generators satisfy the covariance relation $\lambda_s\sigma(a)\lambda_{s}^* = \sigma(\alpha_s(a))$ for all $s \in \Gamma$ and $a \in \mathcal{A}$. Note that $\mathcal{A}\rtimes_{r}\Gamma$ naturally contains the reduced group $C^*$-algebra, $C_r^*(\Gamma) = \overline{\text{Span}\{\lambda_s : s\in\Gamma\}}$, as a subalgebra.

Finally, there exists a canonical $\Gamma$-equivariant conditional expectation $\mathbb{E}:\mathcal{A}\rtimes_{r}\Gamma\to\mathcal{A}$. On the dense subspace of finite sums, it is defined by selecting the coefficient of the identity element:
\[
\mathbb{E}\left(\sigma(a)\lambda_s\right) = 
\begin{cases} 
a & \text{if } s=e, \\
0 & \text{if } s \neq e.
\end{cases}
\] 
\starttocentries
\section{Universal Strongly Proximal Extension for Subgroups}\label{sec:uni}
Having established the rigidity properties of generalized boundaries, we now turn to the question of existence. We consider the category of $\Gamma$-extensions that are minimal as $\Gamma$-spaces but exhibit strong proximality only relative to a subgroup $H$. Using an inverse limit construction (used in \cite{SGlasnerbook}) based on fiber products, we prove the existence of a universal object in this category. In particular, we prove Theorem \ref{thmA}.

	Let $H\le \Gamma$ and let $X$ be a $\Gamma$-space.
Let $F$ be a finite set and for each $\alpha\in F$ let $\pi_\alpha:Y_\alpha\to X$ be a $\Gamma$-extension.
Define the (finite) fibre product
\[
Z_F:=\Bigl\{(x,(y_\alpha)_{\alpha\in F})\in X\times\prod_{\alpha\in F}Y_\alpha:\ 
\pi_\alpha(y_\alpha)=x\ \forall\alpha\in F\Bigr\},
\]
and	 the $\Gamma$-action on $Z_F$ by
\[
t (x,(y_\alpha)_\alpha):=\bigl(tx,(ty_\alpha)_\alpha\bigr),\qquad t\in\Gamma.
\]
Then it is straightforward to check $Z_F$ is a $\Gamma$-space. Define the map
\[
\pi_F:Z_F\to X,\qquad \bigl(x,(y_\alpha)_{\alpha\in F}\bigr)\mapsto x.
\]
Then $\pi_F$ is clearly a  $\Gamma$-factor map. A standard tool for constructing universal objects is to use the fibre product. We first verify that the property of $H$-strong proximality is preserved under taking finite fibre products.
\begin{lemma}\label{lem:finite_fibre_product_HSP}
		Let $H\le \Gamma$ and let $X$ be a $\Gamma$-space.
	Given a finite set $F$,  let $\pi_\alpha:Y_\alpha\to X$ be a $\Gamma$-extension, $\alpha\in F$.
Assume that each $\pi_\alpha$ is an $H$-strongly proximal extension over $X$.
	Then $\pi_F$ is an $H$-strongly proximal extension over $X$.
\end{lemma}
\begin{proof}
	We argue by induction on $|F|$. The case $|F|=1$ is the assumption.
	
	Assume $|F|\ge 2$ and the statement holds for all smaller finite sets.
	Pick $\alpha_0\in F$ and write $F':=F\setminus\{\alpha_0\}$.
	Set
	\[
	Z_{F'}:=\Bigl\{(x,(y_\alpha)_{\alpha\in F'})\in X\times\prod_{\alpha\in F'}Y_\alpha:\ 
	\pi_\alpha(y_\alpha)=x\ \forall\alpha\in F'\Bigr\}.
	\]
	Then $Z_F$ can be identified with the fibre product $Y_{\alpha_0}\times_X Z_{F'}$ via the map
	\[
	(x,(y_\alpha)_{\alpha\in F})\mapsto \bigl(y_{\alpha_0},(x,(y_\alpha)_{\alpha\in F'})\bigr),
	\]
	and under this identification $\pi_F$ is the induced map to $X$. Fix $x\in X$ and let $\nu\in\Prob(Z_F)$ with $\supp(\nu)\subset \pi_F^{-1}(x)$.
	Let $q_1:Z_F\to Y_{\alpha_0}$ and $q_2:Z_F\to Z_{F'}$ be the coordinate projections. Furthermore, note that for every $(x,(y_\alpha)_{\alpha\in F})\in Z_F$ we have
$\pi_{\alpha_0}(y_{\alpha_0})=x$, hence
\[
  \pi_{\alpha_0}\circ q_1(x,(y_\alpha)_{\alpha\in F})
  =\pi_{\alpha_0}(y_{\alpha_0})
  =x
  =\pi_F(x,(y_\alpha)_{\alpha\in F}).
\]
Thus $\pi_{\alpha_0}\circ q_1=\pi_F$ on $Z_F$. Consequently,  $\supp\bigl((q_1)_*\nu\bigr)\subset \pi_{\alpha_0}^{-1}(x)$. 
	By $H$-strong proximality of $\pi_{\alpha_0}$, there exist $y_0\in Y_{\alpha_0}$ and a net $(h_i)\subset H$ such that
	\[
	h_i (q_1)_*\nu \to \delta_{y_0}\quad \text{in }\Prob(Y_{\alpha_0}).
	\]
	Passing to a subnet, we may assume $h_i\nu\to \nu'$ in $\Prob(Z_F)$.
	Then
	\[
	(q_1)_*\nu'=\lim_i (q_1)_*(h_i\nu)=\lim_i h_i (q_1)_*\nu=\delta_{y_0},
	\]
	hence $\nu'$ is supported on $q_1^{-1}(y_0)$.
	
	Since $\pi_F$ is $\Gamma$-equivariant, we have $\supp(h_i\nu)\subset \pi_F^{-1}(h_i x)$ for all $i$.
	As $X$ is compact, after passing to a further subnet, we may assume $h_i x\to x_1\in X$.
	We claim that
	\[
	\supp(\nu')\subset \pi_F^{-1}(x_1).
	\]
	Indeed, suppose for a contradiction that there exists $z\in 	\supp(\nu')\setminus \pi_F^{-1}(x_1)$, and hence there exists an open neighborhood $U$ of $z$ with $\overline{U}\cap \pi_F^{-1}(x_1)=\emptyset$ and $\nu'(U)>0$. Then $x_1\notin \pi_F(\overline{U})$. Since $\pi_F(\overline U)$ is compact, there exists an open neighborhood
	$V$ of $x_1$ with $V\cap \pi_F(\overline U)=\emptyset$.
	For $i$ large enough we have $h_i x\in V$, hence $\pi_F^{-1}(h_i x)\cap \overline U=\emptyset$.
	Therefore $(h_i\nu)(U)=0$ for all large $i$, and by Portmanteau theorem (see e.g. \cite[p. 149 Remarks (3)]{Waltersbook}), one has $\nu'(U)=0$.
	This is a contradiction and proves the claim.

Thus $\nu'$ is supported on
\[
q_1^{-1}(y_0)\cap \pi_F^{-1}(x_1)=\{y_0\}\times \pi_{F'}^{-1}(x_1)
\quad\subset\ Y_{\alpha_0}\times_X Z_{F'}.
\]
Consequently, the measure $\widetilde\nu:=(q_2)_*\nu'\in\Prob(Z_{F'})$ satisfies
$\supp(\widetilde\nu)\subset \pi_{F'}^{-1}(x_1)$.

By the induction hypothesis, $\pi_{F'}:Z_{F'}\to X$ is $H$-strongly proximal over $X$.
Hence there exist $z\in Z_{F'}$ and a net $(k_j)\subset H$ such that
\[
k_j\widetilde\nu \to \delta_{z}\quad \text{in }\Prob(Z_{F'}).
\]
Consider the net $k_j\nu'\in\Prob(Z_F)$. Since $\Prob(Z_F)$ is compact, after passing to a subnet, we may assume
$k_j\nu'\to \nu''$ in $\Prob(Z_F)$. Then
\[
(q_2)_*\nu''=\lim_j (q_2)_*(k_j\nu')=\lim_j k_j(q_2)_*\nu'=\lim_j k_j\widetilde\nu=\delta_z.
\]
On the other hand,
\[
(q_1)_*(k_j\nu')=k_j(q_1)_*\nu'=k_j\delta_{y_0}=\delta_{k_j y_0}.
\]
Passing to a further subnet, we may assume $k_j y_0\to y_1\in Y_{\alpha_0}$, and then
\[
(q_1)_*\nu''=\lim_j (q_1)_*(k_j\nu')=\lim_j \delta_{k_j y_0}=\delta_{y_1}.
\]
Therefore $\nu''=\delta_{(y_1,z)}$ is a Dirac measure on $Z_F$.
Finally, since $\nu' \in \overline{H\nu}$ and $\nu''\in \overline{H\nu'}$, we have $\nu''\in \overline{H\nu}$.
This shows that $\overline{H\nu}$ contains a Dirac measure, completing the induction.
\end{proof}
To capture the universal property, we need to take the limit over all possible extensions. This requires extending the stability result to infinite fibre products.
Similarly,	let $\{ \pi_\alpha:Y_\alpha\to X\}_{\alpha\in A}$ be a family of $\Gamma$-extensions, and
	let
\[
Z:=\Bigl\{(x,(y_\alpha)_{\alpha\in A})\in X\times\prod_{\alpha\in A}Y_\alpha:\ 
\pi_\alpha(y_\alpha)=x \text{ for all }\alpha\in A\Bigr\}.
\]
Then $Z$ is closed in $X\times\prod_{\alpha\in A}Y_\alpha$, hence compact.
Define the $\Gamma$-action on $Z$ by
\[
t (x,(y_\alpha)_\alpha):=\bigl(tx,(ty_\alpha)_\alpha\bigr),\qquad t\in\Gamma,
\]
and define
\[
\pi_Z:Z\to X,\qquad \pi_Z(x,(y_\alpha)_\alpha):=x,
\]
and the coordinate maps
\[
p_\alpha:Z\to Y_\alpha,\qquad p_\alpha(x,(y_\beta)_\beta):=y_\alpha.
\]
Then $\pi_Z$ and all $p_\alpha$ are continuous $\Gamma$-maps, and
\[
\pi_\alpha\circ p_\alpha=\pi_Z\qquad(\forall \alpha\in A).
\]
Moreover, $\pi_Z$ is surjective: given $x\in X$, choose $y_\alpha\in \pi_\alpha^{-1}(x)$ for each $\alpha$
(which is possible since $\pi_\alpha$ is onto), then $(x,(y_\alpha)_\alpha)\in Z$.
\begin{lemma}\label{lem:fibre_product_HSP}
	Let $H\le \Gamma$ and let $X$ be a $\Gamma$-space.
	Let $\{ \pi_\alpha:Y_\alpha\to X\}_{\alpha\in A}$ be a family of $\Gamma$-extensions such that each
	$\pi_\alpha$ is an $H$-strongly proximal extension over $X$.
	Then $\pi_Z:Z\to X$ is an $H$-strongly proximal extension over $X$.
\end{lemma}
\begin{proof}
	Fix $x\in X$ and $\nu\in\Prob(Z)$ with $\supp(\nu)\subset \pi_Z^{-1}(x)$.
	For each finite subset $F\subset A$, let $p_F:Z\to Z_F$ be the projection onto
	\[
	Z_F:=\Bigl\{(x,(y_\alpha)_{\alpha\in F})\in X\times\prod_{\alpha\in F}Y_\alpha:\ 
	\pi_\alpha(y_\alpha)=x\ \forall\alpha\in F\Bigr\},
	\]
	and put $\nu_F:=(p_F)_*\nu\in\Prob(Z_F)$.
	Applying Lemma \ref{lem:finite_fibre_product_HSP} to $\nu_F$, we obtain a Dirac measure $\delta_{z_F}\in \overline{H\nu_F}$.

	Choose a net $(h_i)\subset H$ such that $h_i\nu_F\to \delta_{z_F}$ in $\Prob(Z_F)$.
	Since $\Prob(Z)$ is compact, after passing to a subnet, we may assume $h_i\nu\to \mu_F$ in $\Prob(Z)$.
	Then by continuity of push-forward,
	\[
	(p_F)_*\mu_F=\lim_i (p_F)_*(h_i\nu)=\lim_i h_i (p_F)_*\nu=\lim_i h_i\nu_F=\delta_{z_F}.
	\]
	Define the closed set
	\[
	\mathcal D_F:=\{\mu\in\Prob(Z):\ (p_F)_*\mu \text{ is a Dirac measure on } Z_F\}.
	\]
	Thus $\overline{H\nu}\cap \mathcal D_F\neq\emptyset$ for every finite subset $F\subset A$.

	If $F_1,\dots,F_n$ are finite and $F:=\bigcup_{j=1}^n F_j$, then
	$\mathcal D_F\subset \bigcap_{j=1}^n \mathcal D_{F_j}$.
	Hence the family $\{\overline{H\nu}\cap \mathcal D_F\}_{F\subset A,|F|<\infty}$ has the finite intersection property.
	Since $\overline{H\nu}$ is compact, there exists
	\[
	\mu\in \bigcap_{F\subset A,|F|<\infty}\bigl(\overline{H\nu}\cap \mathcal D_F\bigr).
	\]
	Consequently, we see that \ $\pi_Z$ is $H$-strongly proximal over $X$.
\end{proof}

	Let $X$ be a $\Gamma$-space.
	A \emph{$\Gamma$-ambit} is a pointed $\Gamma$-space $(X,x_0)$ such that the orbit
	$\Gamma x_0:=\{tx_0:t\in\Gamma\}$ is dense in $X$.
We are now in a position to prove the existence of a universal generalized boundary for a subgroup, effectively constructing the relative Furstenberg boundary in this context by using the maximal strongly proximal extension. Recall that $\mathcal C$ is the collection of all $\Gamma$-extensions $\pi:Y\to X$ such that
	\begin{enumerate}
		\item $Y$ is $\Gamma$-minimal;
		\item $\pi$ is an $H$-strongly proximal extension over $X$.
	\end{enumerate}
Moreover, morphisms in $\mathcal C$ are $\Gamma$-factor maps $f:Y\to Y'$
	satisfying $\pi' \circ f=\pi$.	
\begin{proof}[Proof of Theorem \ref{thmA}]
	Since $X$ is $H$-minimal and $H\le \Gamma$, it follows that $X$ is $\Gamma$-minimal. We begin with the following preliminary observation.
\begin{claim}\label{cl3}
There exists a set $\mathscr{A}\subseteq\mathcal{C}$ of representatives such that
every object of $\mathcal{C}$ is isomorphic to some element of $\mathscr{A}$.
\end{claim}
    \begin{proof}[Proof of Claim \ref{cl3}]
        Fix $x_0\in X$. Since $X$ is $\Gamma$-minimal, $(X,x_0)$ is a $\Gamma$-ambit.
	Let $\beta\Gamma$ be the Stone--\v{C}ech compactification of the discrete group $\Gamma$
	with the canonical $\Gamma$-action by left translation, and $e\in\Gamma$, the identity element. Then $(\beta\Gamma,e)$ is the universal $\Gamma$-ambit (see e.g. \cite{Universaladmit}).
	
	For any object $\pi:Y\to X$ in $\mathcal C$, pick $y_0\in\pi^{-1}(x_0)$.
	Since $Y$ is $\Gamma$-minimal, $(Y,y_0)$ is a $\Gamma$-ambit, hence there exists a surjective $\Gamma$-map
	$p:\beta\Gamma\to Y$ with $p(e)=y_0$. Consequently, $Y$ is a $\Gamma$-quotient of $\beta\Gamma$,
	so $Y\cong \beta\Gamma/R$ for some closed $\Gamma$-invariant equivalence relation $R$ on $\beta\Gamma$.
	
	Let $\mathscr R$ be the set of all closed $\Gamma$-invariant equivalence relations on $\beta\Gamma$
	(this is a set since it is a subset of $\mathcal P(\beta\Gamma\times \beta\Gamma)$).
	For each $R\in\mathscr R$, put $Y_R:=\beta\Gamma/R$ with the induced $\Gamma$-action.
	Let $\mathscr M_R$ be the set of all continuous $\Gamma$-maps $\pi:Y_R\to X$
	(it is a set since it is a subset of $X^{Y_R}$).
	Define
	\[
	\mathscr A:=\{(R,\pi): R\in\mathscr R,\ \pi\in\mathscr M_R,\ \pi:Y_R\to X \text{ belongs to }\mathcal C\}.
	\]
	Then $\mathscr A$ is a set, and every object of $\mathcal C$ is $\Gamma$-isomorphic over $X$
	to some $(Y_{(R,\pi)},\pi)$ with $(R,\pi)\in\mathscr A$.
	For simplicity, re-index $\mathscr A$ as a set $A$ and write the corresponding family as
	\[
	\{(\pi_\alpha:Y_\alpha\to X)\}_{\alpha\in A}\subset \mathcal C
	\]
	such that every object of $\mathcal C$ is $\Gamma$-isomorphic over $X$ to some $(Y_\alpha,\pi_\alpha)$. This here finishes Claim \ref{cl3}.
    \end{proof}

	Now, let
	\[
	Z:=\Bigl\{(x,(y_\alpha)_{\alpha\in A})\in X\times\prod_{\alpha\in A}Y_\alpha:\ 
	\pi_\alpha(y_\alpha)=x \text{ for all }\alpha\in A\Bigr\},
	\]
and define the $\Gamma$-action on $Z$ by
	\[
	t (x,(y_\alpha)_\alpha):=\bigl(tx,(ty_\alpha)_\alpha\bigr),\qquad t\in\Gamma.
	\]
Denote
	\[
	\pi_Z:Z\to X,\qquad (x,(y_\alpha)_\alpha)\mapsto x,
	\]
	and the coordinate maps
	\[
	p_\alpha:Z\to Y_\alpha,\qquad (x,(y_\beta)_\beta)\mapsto y_\alpha.
	\]
	Then $\pi_Z$ and all $p_\alpha$ are continuous $\Gamma$-maps, and
	\[
	\pi_\alpha\circ p_\alpha=\pi_Z\qquad(\forall \alpha\in A).
	\]
 By Lemma \ref{lem:fibre_product_HSP}, $\pi_Z$ is $H$-strongly proximal over $X$.

	Let $M\subset Z$ be a nonempty closed $\Gamma$-invariant subset which is $\Gamma$-minimal.
	Set $\pi_M:=\pi_Z|_M:M\to X$. Then $\pi_M$ is a continuous $\Gamma$-map.
	Since $\pi_Z$ is $\Gamma$-equivariant, the image $\pi_Z(M)$ is a nonempty closed $\Gamma$-invariant subset of $X$.
	As $X$ is $\Gamma$-minimal, it follows that $\pi_Z(M)=X$, hence $\pi_M$ is a $\Gamma$-extension.
	Moreover, as $\pi_Z$ is $H$-strongly proximal over $X$, and $M$ is a $\Gamma$-invariant closed subset of $Z$, it follows that $\pi_M$ is $H$-strongly proximal over $X$.
	    Therefore $\pi_M:M\to X$ belongs to $\mathcal C$.
	
We now verify the universality of $\pi_M$.
To this end, we first record the following observation.
Fix $\alpha\in A$. Since $M\subset Z$ is $\Gamma$--invariant,
its image $p_\alpha(M)\subset Y_\alpha$ is a nonempty compact $\Gamma$-invariant subset.
As $Y_\alpha$ is $\Gamma$-minimal, it follows that $p_\alpha(M)=Y_\alpha$.
Consequently, the restriction
\[
p_\alpha|_M:M\to Y_\alpha
\]
is a surjective continuous $\Gamma$-map, and moreover
\[
\pi_\alpha\circ (p_\alpha|_M)=\pi_Z|_M=\pi_M.
\]

Now let $\pi:Y\to X$ be any object in $\mathcal C$.
By the construction of the family $\{(Y_\alpha,\pi_\alpha)\}_{\alpha\in A}$,
there exist $\alpha\in A$ and a $\Gamma$--isomorphism over $X$,
say $\theta:Y_\alpha\to Y$, such that $\pi\circ\theta=\pi_\alpha$.
Define
\[
f:=\theta\circ (p_\alpha|_M):M\to Y.
\]
Then $f$ is a $\Gamma$-factor map, and
\[
\pi\circ f
=\pi\circ\theta\circ (p_\alpha|_M)
=\pi_\alpha\circ (p_\alpha|_M)
=\pi_M.
\]
This proves that $\pi_M:M\to X$ is universal in the category $\mathcal C$.
\end{proof}
\begin{remark}
In the  case where $X$ is a one-point space, the statement reduces to
\cite[Proposition~4.1]{ursu2022relative}.
\end{remark}

\section{Relative Boundary Extension Theorem}\label{sec:rel}
A fundamental problem in the theory of boundaries is to understand how boundary dynamics persist when one passes from a subgroup to the ambient group. 
In this section, we address this extension problem for generalized Furstenberg boundaries. 
Exploiting the rigidity of minimal strongly proximal actions, we construct a canonical and unique extension of the boundary action from a commensurated subgroup to the whole group, assuming the base action is minimal. 
This yields the proof of Theorem~\ref{thmB} stated in the introduction.

\begin{definition}\label{def:commensurated}
Let $\Gamma$ be a group and $H\leq \Gamma$ a subgroup. We say that $H$ is
\emph{commensurated} in $\Gamma$, and write $H\leq_c \Gamma$, if for every
$g\in \Gamma$ the subgroup $H\cap gHg^{-1}$ has finite index in $H$, i.e.
\[
   [H : H\cap gHg^{-1}] < \infty.
\]
Equivalently, $H\cap gHg^{-1}$ has finitely many left (or right)
cosets in $H$.
\end{definition}
The following theorem generalizes \cite[Theorem II.4.4]{SGlasnerbook}, \cite[Lemma 5.2]{breuillard2017c}, and \cite[Theorem 3.1]{li2023c} to the setting of generalized Furstenberg boundaries.
\begin{theorem}\label{thm:extension of action}
\label{thm:relative-Li-Scarparo}
Let  $H\leq_c \Gamma$ be a commensurated
subgroup. Let $X$ be a  minimal $\Gamma$-space and denote by
\[
   \pi_H: \partial_F(H,X)\to X
\]
the generalized Furstenberg boundary of $(H,X)$. Then there exists a unique action of $\Gamma$ on $\partial_F(H,X)$ by
        homeomorphisms,  which we write as
        \[
           \Gamma\overset{T}{\curvearrowright} \partial_F(H,X),\quad (g, y)\mapsto T_g(y),
        \]
        such that:
        \begin{enumerate}
          \item for every $h\in H$ and $y\in \partial_F(H,X)$,
                $T_h(y)=hy$;
          \item the map $\pi_H$ is $\Gamma$-equivariant:
                \[
                    \pi_H\bigl(T_g(y)\bigr) = g\pi_H(y)
                    \quad\text{for all }g\in \Gamma,\ y\in\partial_F(H,X).
                \]
        \end{enumerate}
\end{theorem}
\begin{proof}
    For each $g\in \Gamma$ we set
\[
   H_g := H\cap gHg^{-1},\qquad H_{g^{-1}} := H\cap g^{-1}Hg.
\]
Since $H\leq_c \Gamma$, both $H_g$ and $H_{g^{-1}}$ have finite index in $H$.
Hence, by Theorem~\ref{thm:equaluptofiniteindex}, there exist unique
$(H_g,X)$- and $(H_{g^{-1}},X)$-boundary isomorphisms
\[
   \varphi_g: \partial_F(H,X)\to \partial_F(H_g,X),
\text{ and }
   \varphi_{g^{-1}}: \partial_F(H,X)\to \partial_F(H_{g^{-1}},X),
\]
which are $H_g$ and $H_{g^{-1}}$-equivariant, respectively, and
satisfy
\[
   \pi_{H_g}\circ \varphi_g = \pi_H,\text{ and }
   \pi_{H_{g^{-1}}}\circ\varphi_{g^{-1}} = \pi_H.
\]

Next, consider the conjugation isomorphism
\[
   c_g: H_{g^{-1}}\to H_g,\qquad c_g(h)=ghg^{-1}.
\]
This is a group isomorphism. The actions of $H_{g^{-1}}$ and $H_g$ on
$X$ are intertwined by $c_g$, in the sense that
\[
   c_g(h)x = (ghg^{-1})x,\quad h\in H_{g^{-1}},\ x\in X;
\]
meanwhile, we consider the map 
\[\phi_g:X\to X,\, x\mapsto gx.\]
 Then
by Lemma \ref{lem:rel-iso-boundary}, there exists a unique homeomorphism
\[
   \widetilde c_g: \partial_F(H_{g^{-1}},X)\to \partial_F(H_g,X)
\]
such that:
\begin{enumerate}[(a)]
  \item for all $h\in H_{g^{-1}}$ and $y\in\partial_F(H_{g^{-1}},X)$,
        \[
            \widetilde c_g(hy) = c_g(h) \widetilde c_g(y) = ghg^{-1} \widetilde c_g(y);
        \]
  \item $ \pi_{H_g}\circ \widetilde c_g =g \pi_{H_{g^{-1}}}.$
\end{enumerate}

We now define, for each $g\in \Gamma$,
\[
   T_g := \varphi_g^{-1}\circ \psi_g\circ \varphi_{g^{-1}}
   : \partial_F(H,X)\to \partial_F(H,X).
\]
Since $\psi_g$, $\varphi_{g^{-1}}$ and $\varphi_g$ are homeomorphisms, it follows that $T_g$ is again a homeomorphism. Moreover, for all $g\in \Gamma$,
\[
   \pi_H\circ T_g = \pi_H\circ \varphi_g^{-1}\circ \psi_g\circ \varphi_{g^{-1}}= g\pi_H.
\]
This shows that (ii) holds.

Let $g\in \Gamma$, $h\in H_{g^{-1}}=H\cap g^{-1}Hg$ and
	$x\in\partial_F(H,X)$. Then
	\begin{align*}
		T_g(hx)
		&= \varphi_g^{-1}\Bigl(\psi_g\bigl(\varphi_{g^{-1}}(hx)\bigr)\Bigr)\\
		&= \varphi_g^{-1}\Bigl(\psi_g\bigl(h\varphi_{g^{-1}}(x)\bigr)\Bigr)\\
		&= \varphi_g^{-1}\Bigl(ghg^{-1}\psi_g(\varphi_{g^{-1}}(x))\Bigr)\\
		&= ghg^{-1}\varphi_g^{-1}\bigl(\psi_g(\varphi_{g^{-1}}(x))\bigr)\\
		&= ghg^{-1}T_g(x).
	\end{align*}
	Thus, for all $h\in H\cap g^{-1}Hg$ and all $x\in\partial_F(H,X)$,
	\begin{equation}\label{eq:conjugacy-Tg}
		T_g(hx) = ghg^{-1}T_g(x).
	\end{equation}

   	We now prove $T_{gh}=T_gT_h$ for all $g,h\in \Gamma$.
	Let $g,h\in \Gamma$ and define
	\[
	K := H\cap h^{-1}Hh\cap (gh)^{-1}H(gh).
	\]
	Since $H\leq_c \Gamma$, each intersection $H\cap h^{-1}Hh$ and
	$H\cap (gh)^{-1}H(gh)$ has finite index in $H$, and hence so does $K$.
	In particular $K\leq H$ is a finite-index subgroup.

	First, by~\eqref{eq:conjugacy-Tg} applied to $gh$ instead of $g$ we see
	that for all $k\in H\cap (gh)^{-1}H(gh)$ and $x\in\partial_F(H,X)$,
	\begin{equation}\label{eq:Tgh-K}
	T_{gh}(kx) = (gh)k(gh)^{-1}T_{gh}(x).
\end{equation}
On the other hand, for $k\in K$ and $x\in\partial_F(H,X)$ we
	first use~\eqref{eq:conjugacy-Tg}  to obtain
	\[
	T_h(kx) = hkh^{-1}T_h(x),
	\]
	since $K\subset H\cap h^{-1}Hh$. Note that $hkh^{-1}\in H\cap g^{-1}Hg$, and hence by \eqref{eq:conjugacy-Tg},
		\[
	T_g(hkh^{-1}y) = ghk h^{-1}g^{-1}T_g(y) = (gh)k(gh)^{-1}T_g(y)
	\]
	for all $y\in\partial_F(H,X)$. Taking $y=T_h(x)$ gives
	\begin{equation}\label{eq:TgTh-K}
		T_gT_h(kx)
		= (gh)k(gh)^{-1}T_gT_h(x),
		\quad k\in K,\ x\in\partial_F(H,X).
	\end{equation}
	
	Set
	\[
	\alpha := (T_gT_h)^{-1}\circ T_{gh}
	\colon \partial_F(H,X)\to \partial_F(H,X).
	\]
	Combining~\eqref{eq:Tgh-K} and~\eqref{eq:TgTh-K} we see that for all
	$k\in K$ and $x\in\partial_F(H,X)$,
	\[
	\alpha(kx)
	= (T_gT_h)^{-1}\bigl(T_{gh}(kx)\bigr)
	= (T_gT_h)^{-1}\bigl((gh)k(gh)^{-1}T_{gh}(x)\bigr)
	= k\alpha(x).
	\]
    Thus, $\alpha$ is $K$-equivariant. Moreover, as $k\in K\subset H$, it follows that 
   \begin{equation}\label{eq:20251129}
       \pi_H\circ \alpha=(gh)^{-1}gh\circ\pi_H=\pi_H.
   \end{equation}
    
  Since $[H:K]<\infty$, it follows from Theorem \ref{thm:equaluptofiniteindex} that $\partial_F(H,X)$ and $\partial_F(K,X)$ are isomorphic as $(K,X)$-boundaries. In particular, we may apply Lemma~\ref{lem:relative-Glasner} to the $K$-equivariant map
$
   \alpha.
$
As $\alpha$ satisfies \eqref{eq:20251129}, Lemma~\ref{lem:relative-Glasner} yields
\[
   \alpha = \id_{\partial_F(H,X)}, \text{ that is, } T_{gh} = T_gT_h.
\]

We next check that for any $g\in H$, $T_g=g$. Fix $g\in H$. Then $H\cap g^{-1}Hg=H$. By \eqref{eq:conjugacy-Tg}, one has 
\[T_g(hx)=ghg^{-1}T_g(x),\text{ for any }h\in H\text{ and }x\in \partial_F(H,X).\]
Let 
\[\beta:=g^{-1}T_g:\partial_F(H,X)\to \partial_F(H,X).\]
Then for any $h\in H$,
\[\beta  h=g^{-1}T_gh=g^{-1}ghg^{-1}T_g=h\beta,\]
i.e.\ $\beta$ is $H$-equivariant.
Furthermore, 
\[\pi_H\circ\beta= g^{-1}g\pi_H=\pi_H.\]
Therefore, by Lemma \ref{lem:relative-Glasner}, one has
\[g^{-1}T_g=\id_{\partial_F(H,X)},\text{ that is, }T_g=g.\]

Finally, we prove the uniqueness. If $\{S_g\}_{g\in \Gamma}$ is another action of $\Gamma$ on
	$\partial_F(H,X)$ satisfies (i) and (ii) in the theorem.
	By (ii), we have $\pi_H\circ S_g^{-1}\circ T_g=\pi_H$ for all $g\in \Gamma$, and by (i),  we have for any $h\in H\cap g^{-1}Hg$ and $x\in \partial_F(X,H)$, 
    \[ S_g^{-1}(T_g(hx))= S_g^{-1}(ghg^{-1}T_g(x))=S_g^{-1}(S_{ghg^{-1}}T_g(x))=h S_g^{-1}(T_g(x)).\]
    Using Lemma \ref{lem:relative-Glasner} again, one has  
    \[S_g=T_g,\text{ for all } g\in \Gamma.\]
    Thus, the extension is unique.
\end{proof}
To analyze the dynamical properties of this extended action—particularly its freeness—we need to understand the topology of the fixed-point sets. We recall the following structural result, which is a combination of \cite[Proposition 3.3]{kawabe2017uniformlyrecurrentsubgroupsideal} and Frolík's theorem \cite{Frolik}.
\begin{lemma}\label{lem:fix point is clopen}
    Let $X$ be a minimal $\Gamma$-space and  let $\partial_F(\Gamma,X)$ be the generalized Furstenberg boundary. Then for any $g\in \Gamma$, the set $\mr{Fix}(g)$ of $g$-fixed points in $\partial_F(\Gamma,X)$ is clopen (i.e.\ closed and open).
\end{lemma}
Let $X$ be a $\Gamma$-space and let $U\subset X$ be nonempty.
We denote by
\[
   N_\Gamma(U):=\{g\in\Gamma : gU\cap U\neq\emptyset\}
\]
the \emph{return time set} of $U$ (with respect to $\Gamma$).
 We denote by $L_\Gamma(U)$ the group generated by the return time set $N_\Gamma(U)$. When the acting group is clear from the context, we omit the subscript
and write $N(U)$ and $L(U)$.

Minimality of the base action constrains the recurrence of open sets in the extension. This leads to the following observation regarding the index of return-time groups.
\begin{lemma}\label{lem:return time}
    Let $X$ be a  minimal $\Gamma$-space, and let $U$ be a non-empty open subset of $X$. Then the group $L(U)$
has finite index. 
\end{lemma}
\begin{proof}
Since $X$ is minimal, it follows from  \cite{Furstenbergbook} (see also  \cite[Lemma 2.2]{liu2025independencemeansensitivityminimal}) that $N(U)$ is a {\em syndetic set} in $\Gamma$, i.e. there exists a finite subset $S$ of $\Gamma$ such that 
\[\Gamma=\bigcup_{s\in S}sN(U).\]
  Therefore, 
  \[[\Gamma:L(U)]\le |S|<\infty,\]
  which implies that $L(U)$ has finite index.
\end{proof}
We say that $X$ is \emph{connected} if there do not exist two nonempty open sets
$U,V\subset X$ such that
\[
U\cap V=\varnothing \quad\text{and}\quad X=U\cup V.
\]
When the base space is connected, the topological restrictions become even stronger, forcing the local return-time behavior to globalize. This generalizes \cite[Lemma 5.1]{breuillard2017c}, where this was shown for a singleton $X$.
\begin{lemma}\label{lem:full group under connected}
 Let  $H\leq_c \Gamma$ be a commensurated
subgroup.   Let $X$ be a  minimal connected $\Gamma$-space and denote by
\[
   \pi_H: \partial_F(H,X)\to X
\]
the generalized Furstenberg boundary of $(H,X)$. Given a nonempty open subset $U$ of $\pF(H,X)$. Then $L(U)=\Gamma$.
\end{lemma}
\begin{proof}
    By Lemma \ref{lem:return time}, there exists a finite subset $F$ of $\Gamma$ such that 
    \[\Gamma=\bigsqcup_{s\in F}sL(U),\]
    which implies that $\{sLU\}_{s\in F}$ is a clopen partition of $\pF(H,X)$. 
    
  Since $\pF(H,X)$ is a strongly proximal (hence proximal) extension, if $x,x'$ lie in
distinct atoms of the partition $\{sLU\}_{s\in F}$, then
$\pi_H(x)\neq \pi_H(x')$. By continuity of $\pi_H$, the family
$\{\,s\cdot \pi_H(LU)\,\}_{s\in F}$ is a pairwise disjoint clopen
partition of $\pi_H(\pF(H,X))=X$. As $X$ is connected, this forces
$|F|=1$, and hence $\Gamma=L(U)$.
\end{proof}
The following example shows that the connectedness assumption on $X$ cannot be removed in general.
\begin{example}
\label{exa:odometer-counterexample}
Let $\Gamma=H=\mb Z$ and let $(X,T)$ be the classical $2$-adic odometer, i.e., $X$ is the Cantor group of $2$-adic integers and $T:X\to X$ is the adding-one map $T(x)=x+1$. It is well known that $(X,T)$ is a minimal $\mb Z$-system. Since $X$ is totally disconnected, it is, in particular, not
connected. We can view each element $x\in X$ as $x=x_0x_1\ldots$. Define the clopen subset
\[
   U:=\{x\in X : x_0=0\}.
\]
It is now easy to see that
\[
   T^nU\cap U\neq\emptyset
   \iff n\in 2\mb Z,
\]
because adding an odd integer flips $x_0$,
while adding an even integer preserves it. Hence, $N(U)=\{n\in\mb Z:T^nU\cap U\neq\emptyset\}=2\mb Z$, and therefore 
$L(U)=2\mb Z\subsetneq\mb Z$. Now let $\pi_H:\partial_F(H,X)\to X$ be the generalized Furstenberg boundary of $(H,X)$, and $T$ denotes the $H$-action. Set
\[
   V:=\pi_H^{-1}(U)\subset \partial_F(H,X).
\]
Then $V$ is a nonempty open subset of $\partial_F(H,X)$. If
$T^nV\cap V\ne\emptyset$, we can choose
$y\in V$ with $T^n y\in V$. For $x=\pi_H(y)\in U$ and using the
$H$-equivariance of $\pi_H$, we obtain
\[
   T^n x = \pi_H(T^n y)\in U,
\]
so that $T^nU\cap U\neq\emptyset$. Hence, $n\in N(U)=2\mb Z$. Thus $N(V)\subset N(U)=2\mb Z$, and hence $L(V)=2\mb Z$.
\end{example}
Combining the extension theorem with the topological properties of the fixed-point sets, we arrive at the following dichotomy for the extended action, which relates dynamical freeness to the algebraic structure of the stabilizer.
\begin{theorem}
\label{thm:extensionaction}
     Let $X$ be a $\Gamma$-space and $H\leq_c\Gamma$ a commensurated subgroup such that
\begin{enumerate}[(1)]
  \item  $X$ is $H$-minimal; 
  \item the action $H\curvearrowright\partial_F(H,X)$ is free.
\end{enumerate}
Let $\{T_g\}_{g\in G}$  denote the homeomorphisms given by Theorem~\ref{thm:extension of action}.
Then for every $g\in\Gamma$ exactly one of the following alternatives holds:
\begin{enumerate}[(i)]
  \item $T_g$ is free, i.e.\ $\Fix(T_g)=\emptyset$;
  \item $g\in C_\Gamma\big(L_{H \cap g^{-1}Hg}(\Fix(g))\big)$.
\end{enumerate}
Moreover:
\begin{itemize}
  \item if $X$ is connected, then $L_{H \cap g^{-1}Hg}(\Fix(g))=H\cap g^{-1}Hg$;
  \item if $g$ acts trivially on $X$, then $g\in C_\Gamma\big(L_{H \cap g^{-1}Hg}(\Fix(g))\big)$ implies $T_g=\id_{\pF(H,X)}$.
\end{itemize}
\end{theorem}
\begin{proof}
 We suppose that the action $T_g$ is not free, and prove  $g\in C_\Gamma\big(L(\Fix(g))\big)$. Since $H\curvearrowright\partial_F(H,X)$ is a strongly proximal and minimal extension over $X$ and $H\cap g^{-1}Hg$ has finite index in $H$, it follows from Lemma~\ref{lem:minimal preserved by normal subgroups} that $H\cap g^{-1}Hg\curvearrowright\partial_F(H,X)$ is minimal. By Lemma \ref{lem:fix point is clopen}, $\Fix(g)$ is a non-empty clopen subset.  Applying Lemma \ref{lem:return time} to $U=\Fix(g)$ and $\Gamma=H\cap g^{-1}Hg$, we see that 
 \begin{equation*}\label{eq:2025122}
     [H \cap g^{-1}Hg: L_{H \cap g^{-1}Hg}(\Fix(g))]<\infty.
 \end{equation*}
For any $k\in N_{H \cap g^{-1}Hg}(\Fix(g))$ and $y\in k\Fix(g)\cap \Fix(g)$, one has 
\[T_g(y)=y,\text{ and }T_g(k^{-1}y)=k^{-1}y,\]
which combining with equation~\eqref{eq:conjugacy-Tg}, implies that 
\[k^{-1}y= T_g(k^{-1}y)=gk^{-1}g^{-1}T_g(y)=gk^{-1}g^{-1}y.\]
Note that $k,gk^{-1}g^{-1}\in H$. Since $H\curvearrowright\partial_F(H,X)$ is free, we get that $kg=gk$. This implies that $g\in C_\Gamma(L_{H \cap g^{-1}Hg}(\Fix(g)))$.

If $X$ is connected, then by Lemma \ref{lem:full group under connected}, we immediately obtain that $$L_{H \cap g^{-1}Hg}(\Fix(g))=H \cap g^{-1}Hg.$$
Finally,  we assume that $g=\id_X$ and $g\in C_\Gamma(L(\Fix(g)))$ , and prove that $T_g=\id_{\pF(H,X)}$. Using \eqref{eq:conjugacy-Tg} again, one has that for any $l\in L_{H \cap g^{-1}Hg}(\Fix(g))$ and $y\in \pF(H,X)$, 
\[T_g(ly)=glg^{-1}T_g(y)=lT_g(y).\]
Thus, $T_g$ is $L_{H \cap g^{-1}Hg}(\Fix(g))$-equivariant.

Since $[H \cap g^{-1}Hg:L_{H \cap g^{-1}Hg}(\Fix(g))]<\infty$ and $[H:H \cap g^{-1}Hg]<\infty$, it follows that $[H:L_{H \cap g^{-1}Hg}(\Fix(g))]<\infty$. By Lemma \ref{lem:key lemma}, one has $\pF(H,X)$ is a $(L_{H \cap g^{-1}Hg}(\Fix(g)),X)$-boundary. Since $T_g$ is an $L_{H \cap g^{-1}Hg}(\Fix(g))$-equivariant continuous self-map of $\pF(H,X)$ with $\pi_H\circ T_g=g\circ \pi_H=\pi_H$, using Lemma \ref{lem:relative-Glasner}, $T_g=\id_{\pF(H,X)}$.
\end{proof}
We conclude this section by examining several concrete classes of groups and actions. We analyze semidirect products and groups acting on trees to identify specific instances where the conditions for boundary extension and trivial centralizers are met.

Let us begin with some notation.
Let $H$ and $N$ be groups, and let $\phi\colon H\to \Aut(N)$ be a homomorphism.
Define a binary operation on $N\times H$ by
\[
(n,h)\cdot (n',h') := \bigl(n\,\phi(h)(n'),\, hh'\bigr),
\qquad (n,h),(n',h')\in N\times H.
\]
With this operation, $N\times H$ becomes a group, called the \emph{semidirect
product} of $N$ by $H$ (with respect to $\phi$), and denoted by
\[
\Gamma := N\rtimes_{\phi} H.
\]
We identify $N$ with the subgroup $N\times\{e_H\}\le \Gamma$ via the embedding
\[
\iota\colon N\to \Gamma,\qquad n\mapsto (n,e_H).
\]
It is straightforward to verify that $N$ is a normal subgroup of
$\Gamma$.

For $n\in N$, let $\Ad_n\in \Aut(N)$ be the inner automorphism given by
\[
\Ad_n(x):=nxn^{-1},\qquad x\in N.
\]
The set of all such inner automorphisms forms a subgroup of $\Aut(N)$, denoted
by $\Inn(N)$. The assignment
\[
\Ad\colon N\to \Inn(N),\qquad n\mapsto \Ad_n,
\]
is a group homomorphism, whose kernel is the center $Z(N):=C_N(N)$.

\begin{proposition}\label{prop:trivial-centralizer-semidirect}
Let $N$ and $H$ be groups, and let $\phi:H\to \Aut(N)$ be a homomorphism.
Assume that
\[
Z(N)=\{e_N\}
\quad\text{and}\quad
\phi^{-1}(\Inn(N))=\{e_H\}.
\]
Let $\Gamma:=N\rtimes_{\phi} H$ be the semidirect product. Then
\[
C_{\Gamma}(N)=\{e_\Gamma\}.
\]
Moreover, if $N$ and $H$ are amenable, then so is $\Gamma$.
\end{proposition}
\begin{proof}
Let $(n,h)\in C_\Gamma(N)$. Then for every $x\in N$ we have
\[
(n,h)(x,e_\Gamma)=(x,e_\Gamma)(n,h).
\]
Compute both sides:
\[
(n,h)(x,e_\Gamma)=\bigl(n\phi(h)(x),\,h\bigr),
\qquad
(x,e_\Gamma)(n,h)=\bigl(xn,\,h\bigr).
\]
Hence
\[
n\phi(h)(x)=xn \quad \forall\,x\in N.
\]
Equivalently,
\[
\phi(h)(x)=n^{-1}xn \quad \forall\,x\in N,
\]
so $\phi(h)=\Ad_{n^{-1}}\in\Inn(N)$.

By the assumption $\phi^{-1}(\Inn(N))=\{e_H\}$, this forces $h=e_H$.
Plugging $h=e_H$ back into $n\phi(h)(x)=xn$ gives $nx=xn$ for all $x\in N$, i.e. $n\in C_N(N)$.
Since $Z(N)=\{e_N\}$, we get $n=e_N$. Therefore $(n,h)=(e_N,e_H)$, so $C_\Gamma(N)=\{e_\Gamma\}$.
\end{proof}
We now provide some concrete examples which guarantee the conditions in Proposition \ref{prop:trivial-centralizer-semidirect}.
\begin{example}\label{exa:Sn-An-C2}[Symmetric groups]
Let $n\ge 5$, let $N=A_n$ be the alternating group, and let $H=C_2=\langle \tau\rangle$ be a 2-cyclic group.
Fix an odd involution $s\in S_n\setminus A_n$ (e.g.\ $s=(1\,2)$) and define
\[
\phi:C_2 \to \Aut(A_n),\qquad \phi(\tau)(x):=sxs^{-1}.
\]
Then $\phi(\tau)^2=c_{s^2}|_{A_n}=\id_{A_n}$, hence $\phi$ is a homomorphism. 
Let $\Gamma:=N\rtimes_\phi H$. Moreover, $\Gamma$ is the symmetric group $S_n$ on $n$ letters.

We now verify  the assumptions in Proposition~\ref{prop:trivial-centralizer-semidirect} to obtain that $C_\Gamma(N)=\{e_\Gamma\}$. Indeed, since $A_n$ is a non-abelian simple group, its center is trivial. Thus,
it suffices to show that $\phi(\tau)\notin \Inn(A_n)$.
Choose $x\in A_n$ as follows:
\[
x=
\begin{cases}
(1\,2\,\dots\,n), & \text{if $n$ is odd},\\[2mm]
(1\,2\,\dots\,n-1), & \text{if $n$ is even}.
\end{cases}
\]
Note that
$sxs^{-1}$ is $A_n$-conjugate to $x$, as $s$ is an odd permutation.

If $\phi(\tau)\in \Inn(A_n)$, then there exists $a\in A_n$ such that
$sxs^{-1}=\phi(\tau)(x)=axa^{-1}$, and hence $sxs^{-1}$ is not $A_n$-conjugate to $x$. This is a contradiction.

Consequently, Proposition~\ref{prop:trivial-centralizer-semidirect} yields
\[
C_\Gamma(N)=\{e_\Gamma\}.
\]
\end{example}

Let $K$ be a group and $I$ a set. The \emph{restricted direct sum}
\[
\bigoplus_{i\in I} K
\]
is the subgroup of $\prod_{i\in I} K$ consisting of all families
$x=(x_i)_{i\in I}$ with \emph{finite support}, i.e.\ the set
$\{i\in I:\ x_i\neq e_K\}$ is finite.
For each $i\in I$ we denote by $K_i\le \bigoplus_{j\in I} K$ the copy of $K$
supported at $i$.

The following lemma is not difficult to prove.
\begin{lemma}\label{lem:inner-preserves-coordinates}
Let $K$ be a group and $I$ a set. Put $N:=\bigoplus_{i\in I} K$.
Then for every $n\in N$ and every $i\in I$, the inner automorphism
$\Ad_n\in\Inn(N)$ satisfies
\[
\Ad_n(K_i)=K_i.
\]
\end{lemma}
\begin{proof}
Fix $n=(n_j)_{j\in I}\in N$ and $x\in K_i$. In $N$ we have
\[
\Ad_n(x)=nxn^{-1}.
\]
Since $x$ is supported only at $i$, the product $nxn^{-1}$ is also supported
only at $i$, and its $i$-th coordinate equals $n_i x_i n_i^{-1}\in K$.
Hence $\Ad_n(x)\in K_i$. So $\Ad_n(K_i)\subseteq K_i$, and equality follows
because $\Ad_n$ is bijective.
\end{proof}
\begin{example}[Shift semidirect products]\label{exa:shift-semidir}
Let $K$ be an amenable group with $Z(K)=\{e_K\}$ (e.g.\ finite group $K=A_m$, $m\ge5$,
or infinite dihedral group $D_\infty$). Set
\[
N:=\bigoplus_{i\in\mb Z} K.
\]
Let $H=\mb Z$, and define the shift
\[
\phi:H\to\Aut(N),\qquad
\phi(1)\big((x_i)_{i\in\mb Z}\big):=(x_{i-1})_{i\in\mb Z}.
\]
Let $\Gamma:=N\rtimes_\phi H$.

We now check the condition in Proposition \ref{prop:trivial-centralizer-semidirect}.
Since $Z(N)=\bigoplus_{i\in \mb Z} Z(K)$ and $Z(K)=\{e_K\}$, it follows that 
\[Z(N)=\{e_N\}.\]

Fix $m\in\mb Z$. Then $\phi(m)$ sends $K_0$ onto $K_m$.
If $m\neq 0$, then $K_m\neq K_0$, so $\phi(t^m)(K_0)\neq K_0$.
On the other hand, by Lemma~\ref{lem:inner-preserves-coordinates},
every inner automorphism of $N$ preserves $K_0$.
Hence $\phi(t^m)\notin\Inn(N)$ for all $m\neq 0$, i.e.\ $\phi^{-1}(\Inn(N))=\{0\}$.

Consequently, Proposition~\ref{prop:trivial-centralizer-semidirect} yields
\[
C_\Gamma(N)=\{e_\Gamma\}.
\]
\end{example}
\begin{remark}\label{rem:generalization}
In Example~\ref{exa:shift-semidir}, one may replace $\mb Z$ by $\mb Z^d$ and $H=\mb Z^d$,
letting $H$ act by translations on $\mb Z^d$ and defining $\phi$ by the corresponding shifts.
All verifications are identical.
\end{remark}
We next provide an example that is not amenable and not $C^*$-simple.
\begin{example}\label{exa:not-Cstar-simple-trivial-centralizer}
Let $K$ by a
nontrivial amenable group $K$ with $Z(K)=\{e\}$ and let $H$ be a nontrivial non-amenable group (e.g. the free group $\mb F_n$).
Put
\[
N:=\bigoplus_{g\in H} K,
\]
and define $\phi:H\to \Aut(N)$ by the left-translation shift
\[
\phi(h)\bigl((x_g)_{g\in H}\bigr):=(x_{h^{-1}g})_{g\in H}\qquad(h\in H).
\]
Let $\Gamma:=N\rtimes_\phi H$.
Then the normal subgroup $N$ satisfies
$
C_\Gamma(N)=\{e_\Gamma\},
$
and $\Gamma$ is not amenable and not $C^*$-simple.
\end{example}

\begin{proof}
For $h_1,h_2\in H$ and $x=(x_g)_{g\in H}\in N$ we have
\[
\phi(h_1)\phi(h_2)(x)
=\phi(h_1)\bigl((x_{h_2^{-1}g})_{g}\bigr)
=(x_{h_2^{-1}h_1^{-1}g})_{g}
=(x_{(h_1h_2)^{-1}g})_{g}
=\phi(h_1h_2)(x).
\]
Hence $\phi:H\to\Aut(N)$ is a homomorphism.

Let $z=(z_g)_{g\in H}\in Z(N)$. Fix $g_0\in H$ and $a\in K$, and let $u\in N$ be supported
only at $g_0$ with $u_{g_0}=a$. Since $zu=uz$, looking at the $g_0$-coordinate gives
$z_{g_0}a=az_{g_0}$ for all $a\in K$, hence $z_{g_0}\in Z(K)=\{e_K\}$.
Thus $z_{g_0}=e_K$ for all $g_0\in H$, so $z=e_N$ and $Z(N)=\{e_N\}$.
By an argument similar to that of Example \ref{exa:shift-semidir}, one has $\phi^{-1}(\Inn(N))=\{e_H\}$. 

Therefore Proposition~\ref{prop:trivial-centralizer-semidirect} applies and yields
$C_\Gamma(N)=\{e_\Gamma\}$.

Since $K$ is amenable, it follows that $N$ is amenable. As $N\lhd \Gamma$ and $N\neq \{e\}$, the amenable radical
 (the largest amenable normal subgroup of $\Gamma$) is nontrivial.
It is a standard fact that $C^*$-simple groups have trivial amenable radical; hence
$\Gamma$ cannot be $C^*$-simple. On the other hand, since $H$ is not amenable, we see that $\Gamma$ is not amenable.
\end{proof}

We now present an example, inspired by \cite{LeBoudec2016,LeBoudec2017}, of a group $\Gamma$
which is not $C^*$-simple and admits a proper normal subgroup $N\lhd \Gamma$
with a trivial centralizer.

\begin{example}\label{ex:LB-proper-normal-trivial-centralizer}
We briefly recall the family of groups introduced by Le~Boudec \cite{LeBoudec2017}.
Let $\Omega$ be a (finite or countable) set with $|\Omega|=d\ge 3$, and let $T_\Omega$
be the $d$-regular tree equipped with a legal edge-colouring by $\Omega$.
Let $F\le F'\le S_\Omega$ be permutation groups.

Following \cite{LeBoudec2016,LeBoudec2017}, define $G(F,F')\le \Aut(T_\Omega)$ to be the group
of all tree automorphisms whose \emph{local action} belongs to $F'$ at every vertex,
and belongs to $F$ at all but finitely many vertices
(see \cite[Section 5]{LeBoudec2017} for a formal definition).
Assume that $F$ acts freely on $\Omega$, preserves the $F'$-orbits, and that all point stabilizers
in $F'$ are amenable. Then $G(F,F')$ is not $C^*$-simple and its amenable radical is trivial
(by \cite[Theorem~C]{LeBoudec2017}).

Assume moreover that $F$ acts transitively on $\Omega$. Let $G(F,F')^{*}\le G(F,F')$ be the
index-two subgroup consisting of automorphisms preserving the natural bipartition (types)
of the vertex set of $T_\Omega$ (cf.\ \cite[after Theorem~C]{LeBoudec2017}).
For suitable choices of $(F,F')$, the group $G(F,F')^{*}$ is (non-abelian) simple
(see \cite[Theorem~4.13 and Corollary~4.14]{LeBoudec2016}).
For instance, if $d\ge 5$ is odd, $F'=A_d$, and $F$ is generated by a $d$-cycle (so that
$F$ is simply transitive), then $G(F,F')^{*}$ is simple \cite[Example~4.15]{LeBoudec2016}.
In particular, in this situation
\[
Z\bigl(G(F,F')^{*}\bigr)=\{e_{G(F,F')^{*}}\}.
\]

Fix such a pair $(F,F')$ and put
\[
\Gamma:=G(F,F'),\qquad N:=G(F,F')^{*}.
\]
Then $N\lhd \Gamma$ is a proper normal subgroup (indeed $[\Gamma:N]=2$), $\Gamma$ is not $C^*$-simple,
and moreover $C_\Gamma(N)=\{e_N\}$ by the following lemma.
\end{example}

\begin{lemma}\label{lem:index2-centralizer-trivial}
Let $\Gamma$ be a group and let $N\le \Gamma$ satisfy $[\Gamma:N]=2$.
Assume that $Z(N)=\{e_N\}$ and that the amenable radical of $\Gamma$ is trivial.
Then
\[
C_\Gamma(N)=\{e_\Gamma\}.
\]
\end{lemma}
\begin{proof}
Put $C:=C_\Gamma(N)$. Since $[\Gamma:N]=2$, we have $N\lhd \Gamma$, hence $C\lhd \Gamma$.
Moreover,
\[
C\cap N =Z(N)=\{e_N\}.
\]
Let $\pi:\Gamma\to \Gamma/N$ be the quotient map. Then
\[
\ker(\pi|_C)=C\cap N=\{e_C\},
\]
so $\pi|_C$ is injective. Since $|\Gamma/N|=2$, it follows that $|C|\le 2$, and hence $C$ is amenable.
Since the amenable radical of $\Gamma$ is trivial, it follows that $C=\{e_\Gamma\}$.
\end{proof}

Let us adapt our notation to reflect the relative nature of these boundaries better. We recall $RB_X(\Gamma, H)$  is the universal object  obtained in Theorem \ref{thmA}. Explicitly, $RB_X(\Gamma, H)$ is the universal $\Gamma$-minimal extension of $X$ which is $H$-strongly proximal with respect to $X$. We now address the relationship between this universal space and the generalized Furstenberg boundary of the subgroup $H$. In the case where $H$ is commensurated in $\Gamma$, we show that they coincide.

\begin{theorem}\label{thm:same}
    Let $H \le_c \Gamma$ be a commensurated subgroup. Suppose that $X$ is an $H$-minimal $\Gamma$-space.  Then there is a $\Gamma$-equivariant homeomorphism
    \[
    RB_X(\Gamma, H) \cong \partial_F(H, X),
    \]
    where $\partial_F(H, X)$ is equipped with the unique $\Gamma$-action constructed in Theorem \ref{thm:relative-Li-Scarparo}.
\end{theorem}

\begin{proof}
    By Theorem \ref{thm:relative-Li-Scarparo}, we know that the generalized Furstenberg boundary $\partial_F(H, X)$ admits a continuous $\Gamma$-action extending the $H$-action. Furthermore, as an $H$-space, $\partial_F(H, X)$ is $H$-minimal and is a relative $H$-boundary over $X$. Since $H \le \Gamma$, the $H$-minimality implies $\Gamma$-minimality. Thus, the space $\partial_F(H, X)$ belongs to the class $\mathcal{C}$ defined in Theorem \ref{thmA} (being $\Gamma$-minimal and $H$-strongly proximal). By the universality of $RB_X(\Gamma, H)$, there exists a unique $\Gamma$-equivariant factor map
    \[
    \phi: RB_X(\Gamma, H) \to \partial_F(H, X).
    \]
    
    On the other hand, consider the space $RB_X(\Gamma, H)$ as an $H$-space. By definition, the map $\pi_{RB}:RB_X(\Gamma, H) \to X$ is an $H$-strongly proximal extension. Moreover, since $RB_X(\Gamma, H)$ is $\Gamma$-minimal and $\pi_{RB}$ is surjective onto the $H$-minimal space $X$, $RB_X(\Gamma, H)$ contains a unique $H$-minimal set (see Lemma \ref{lem:unique of minimal set}). We observe that $RB_X(\Gamma, H)$ is a boundary in the sense of the category of $H$-extensions. By the universal property of the generalized Furstenberg boundary $\partial_F(H, X)$ (Definition \ref{def:gen-F-boundary}), there exists an $H$-equivariant map
    \[
    \psi: \partial_F(H, X) \to RB_X(\Gamma, H).
    \]
    Consider the composition $\Phi:=\phi \circ \psi: \partial_F(H, X) \to \partial_F(H, X)$. Since $\phi$ is $\Gamma$-equivariant, it is $H$-equivariant. Thus $\Phi$ is an $H$-equivariant self-map of $\partial_F(H, X)$ commuting with the factor map to $X$. Using Lemma~\ref{lem:relative-Glasner}, $\Phi$ must be the identity map on $\partial_F(H, X)$. Consequently, $\psi$ is an embedding and $\phi$ is a left inverse to $\psi$. Since $\partial_F(H,X)$ is minimal, it follows that $\phi$ is a homeomorphism. Thus, the spaces are $\Gamma$-equivariantly homeomorphic.
\end{proof}
\begin{remark}
We would like to point out that, a priori, it is not clear whether
$H\curvearrowright RB_X(\Gamma,H)$ is minimal. However, in the case when
$H\le_c\Gamma$, it follows from the above theorem that this is the case.
\end{remark}
\section{\texorpdfstring{$X$}{}-plump}\label{sec:application}
The dynamical behavior of the generalized Furstenberg boundary is deeply intertwined with the ideal structure of the associated reduced crossed products. In this section, we translate our dynamical results into the language of operator algebras by introducing the notion of an $X$-plump subgroup—a relative analogue of the generalized Powers averaging property~\cite{amrutam2022generalized}. We demonstrate that $X$-plumpness is a powerful obstruction to the existence of non-trivial ideals, leading to simplicity results for intermediate $C^*$-algebras and rigidity theorems that characterize when such algebras must themselves be crossed products.

\begin{definition}[$X$-plump]
\label{def:genplump}
Let $X$ be a $\Gamma$-space. We say
that $H<\Gamma$ is $X$-plump if given $\varepsilon>0$ and $a\in C(X)\rtimes_r\Gamma$, we can find
$s_{1},\dots,s_{m}\in H$
and $g_1,\dots,g_m\in C(X)^+$ with $\sum_{j=1}^m g_j^2 = 1$ such that
\[
\left\| \sum_{j=1}^m g_j\, \lambda(s_j)a\, \lambda(s_j)^*\, g_j -\mathbb{E}(a)\right\| < \varepsilon.
\]
\end{definition}

The dynamical properties of the boundary have direct consequences for the analytic structure of the group $C^*$-algebras. The following theorem establishes a connection between the boundary action and the notion of $X$-plumpness, thereby bridging the gap between dynamics and operator algebras. 
\begin{theorem}
\label{thm:condgenplump}
Let $X$ be a $\Gamma$-space and $H\le\Gamma$. Suppose there is
a topologically free action of $\Gamma$ on a $(\Gamma,X)$-boundary $\partial(\Gamma,X)$
such that 
\begin{enumerate}
    \item $H\curvearrowright X$ is minimal.
    \item the action $\Gamma \curvearrowright \partial(\Gamma,X)$ restricted
to $H$ is strongly proximal and minimal with respect to $X$.
\end{enumerate}
Then
$H$ is $X$-plump.
\end{theorem}
\begin{proof}
We argue similarly as in \cite{amrutam2022generalized} and adapt it to our purposes. Given any state $\phi \in S(C(X) \rtimes_r \Gamma)$, we claim that
    \[ \left\{\nu \circ \bE: \nu \in \text{Prob}(X)\right\} \subseteq \overline{\left\{\phi \mu:\mu \in P_f(H,C(X))\right\}}^{\text{weak}^*}. \] 
    Recall that $P_f(H,C(X))$ is the collection of all \emph{generalized $(C(X), H,\alpha)$-probability measures}, which is a finitely supported function $f \colon H \to C(X)_+$ which satisfies $\sum_{s \in H} f(s)^2 = 1$ (where the sum converges uniformly). To do so, it is enough to show that
    $$ \delta_x \circ \bE \in \overline{ \{ \phi\mu : \mu \in P_f(H, C(X)) \} } \quad \forall x \in X. $$
Given a state $\phi \in S(C(X) \rtimes_r \Gamma)$, extend it to a state $\tilde{\phi}$ on $C(\partial(\Gamma, X) \rtimes_r \Gamma)$ and let
$$ \eta = \tilde{\phi}|_{C(\partial(\Gamma, X))} $$
Then since $H \curvearrowright X$ is minimal, using \cite[Lemma~3.6]{amrutam2022generalized}, $\exists \mu \in \text{Prob}(H, C(X))$ such that $\eta \mu \to \delta_{x'}$ for some $x' \in X$.
Upon passing to a subnet, assume $\eta \mu \to \tilde{\eta}$ (where $\pi_* \tilde{\eta} = \delta_x$). Since $H \curvearrowright \partial(\Gamma, X)$ is strongly proximal and minimal, $\exists$ a net $\{t_i\} \subseteq H$ such that $t_i \tilde{\eta} \to \delta_y$, where $y \in \partial(\Gamma, X)$ is a free point in the sense that $ty \neq y$ for all non-identity elements $t \in \Gamma$. Then, by arguing similarly as in \cite[Proposition~3.7]{amrutam2022generalized}, we see that
$$ \delta_{\pi(y)} \circ \bE \in \overline{ \{ \phi\mu : \mu \in P_f(H, C(X)) \} }^{w^*} $$
However, since $H \curvearrowright X$ is minimal, and the latter set is convex, the claim follows.
    Consequently, arguing as in the implication of $(1)\implies (4)$ of \cite[Theorem~1.2]{amrutam2022generalized}, we see that if $a \in C(X) \rtimes_r \Gamma$ and $x \in X$, then
    \[ \bE(a)(x) \in \overline{\setbuilder{\mu a}{\mu \in P_f(H,C(X))}}^{\|\cdot\|}. \]  
We can now approximate $\bE(a)$ by $C(X)$-convex combinations of $\bE(a)(x)$, where $x \in X$, arguing similarly as in the implication $(4)\implies (3)$ of \cite[Theorem~1.2]{amrutam2022generalized}. The claim follows By \cite[Remark~2.10]{amrutam2022generalized}, we have that $\bE(a) \in \overline{\setbuilder{\mu a}{\mu \in P_f(H,C(X))}}^{\|\cdot\|}$, as this set is closed under $C(X)$-convex combinations.    
\end{proof}
If $X$ is connected, it becomes much easier, especially when $H$ is a commensurated subgroup.
\begin{corollary}\label{cor:connected} Let $X$ be a connected $\Gamma$-space and $H \leq_c \Gamma$. Suppose that $C_{\Gamma}(H) = \{e\}$. Assume that $H \curvearrowright X$ is minimal and $C(X) \rtimes_r H$ is simple. Then $H$ is $X$-Plump.
\end{corollary}
\begin{proof} Since $C(X) \rtimes_r H$ is simple, $H \curvearrowright\partial_F(H, X)$ is free. Moreover, using Theorem~\ref{thm:extension of action}, $H \curvearrowright\partial_F(H, X)$ extends to an action $\Gamma \curvearrowright \partial_F(H, X)$ such that the action $\Gamma \curvearrowright \partial_F(H, X)$ is free. The claim now follows from Theorem~\ref{thm:condgenplump}.
\end{proof}
The structural properties of intermediate $C^*$-algebras associated with $\Gamma$ were previously analyzed in our earlier work \cite{amrutam2021intermediate} (see also \cite{amrutam2023powers}). Here, we show that $\Gamma$ satisfies the plumpness condition relative to its boundary action, thereby generalizing the simplicity results for intermediate subalgebras. Note that if $X$ is trivial, then we get that $H$ is Plump in $\Gamma$, which forces $\Gamma$ to be $C^*$-simple.  

To apply averaging arguments effectively, we need a preliminary estimate on convolutions of non-negative functions over minimal spaces.
\begin{lemma}\label{lem:connonnegfun}
Let $X$ be a minimal $\Gamma$-space and let $\mu\in\Prob(\Gamma)$ be a
probability measure on $\Gamma$ with full support. If
$f\in C(X)$ is nonnegative and not identically zero, then there exists
$\delta>0$ such that
\[
   \mu*f \ge \delta \quad\text{on }X.
\]
\end{lemma}
\begin{proof}
Since $f\ge 0$ and $\mu$ has full support, the convolution
\[
   (\mu*f)(x)\;=\;\sum_{s\in\Gamma}\mu(s)\,f(s^{-1}x),\qquad x\in X,
\]
is a continuous nonnegative function on $X$. As $X$ is compact, it suffices
to show that $(\mu*f)(x)>0$ for every $x\in X$; then the minimum of
$\mu*f$ on $X$ will be a strictly positive constant~$\delta$.

Assume by contradiction that there exists $x\in X$ with $(\mu*f)(x)=0$. Then
\[
   0 = (\mu*f)(x) = \sum_{s\in\Gamma}\mu(s)\,f(s^{-1}x).
\]
Each term $\mu(s)f(s^{-1}x)$ is nonnegative, and $\mu(s)>0$ for all
$s\in\Gamma$ by full support, so we must have
\[
   f(s^{-1}x)=0 \quad\text{for all }s\in\Gamma.
\]
 Since the action
$\Gamma\curvearrowright X$ is minimal, the orbit $\Gamma x$ is dense in $X$,
and by continuity of $f$ it follows that $f\equiv 0$ on $X$, contradicting
our assumption. Hence $(\mu*f)(x)>0$ for all $x\in X$.
\end{proof}

We can now apply the machinery of $X$-plumpness to study the ideal structure of intermediate $C^*$-algebras. The following result establishes simplicity for a wide class of such inclusions (in particular, we prove Theorem \ref{thmC}). This generalizes \cite[Theorem~1.5]{amrutam2022generalized} and also gives new examples of $C^*$-irreducible inclusions in the sense of \cite{rordam2023irreducible} (also see~\cite{amrutam2020simplicity,kwasniewski2022aperiodicity}).
\begin{theorem}\label{thm:intsim}
Let $\pi: X \to Y$ be a continuous $\Gamma$-equivariant surjection, where
$X$ and $Y$ are compact Hausdorff $\Gamma$-spaces. Suppose there exists a
subgroup $H\le \Gamma$ such that
 $H\curvearrowright X$ is minimal, and
 $H$ is $Y$-plump.
Then every intermediate $C^*$-algebra
\[
   C(Y)\rtimes_r H\ \subseteq\ \mathcal{A}\ \subseteq\ C(X)\rtimes_r\Gamma
\]
is simple.
\end{theorem}
\begin{proof}
Let $\mathcal{A}$ be such an intermediate $C^*$-subalgebra, and let
$I\subseteq\mathcal{A}$ be a nonzero ideal. We aim to show $I=\mathcal{A}$.

Let $0\ne b\in I$ be arbitrary and set $a=b^*b\in I$. Then $a\ge 0$ and
$a\ne 0$. Let $\bE:C(X)\rtimes_r\Gamma\to C(X)$ denote the canonical
conditional expectation. Since $\bE$ is faithful, $\bE(a)\in C(X)^+$ is
nonzero.

Choose a probability measure $\mu\in\Prob(H)$ with full support
$\supp(\mu)=H$. Applying Lemma~\ref{lem:connonnegfun} with the $H$-action
on $X$ (which is minimal by assumption) and the nonzero function
$f=\bE(a)\ge 0$, we obtain a constant $\delta>0$ such that
\[
   (\mu*\bE(a))(x)\ >\ \delta\qquad\text{for all }x\in X,
\]
where
\[
   (\mu*f)(x) := \sum_{s\in H} \mu(s)\, f(s^{-1}x),\quad f\in C(X).
\]

For this fixed $\mu$, define
\[
   \mu*a := \sum_{s\in H} \mu(s)\,\lambda(s)\,a\,\lambda(s)^*
   \ \in C(X)\rtimes_r\Gamma,
\]
where $\lambda(s)$ denotes the canonical unitary implementing the action.
Note that for each $s\in H$, we have $\lambda(s)\in C(Y)\rtimes_r H
\subseteq\mathcal{A}$, and since $a\in I$ and $I$ is a (two-sided) ideal
of $\mathcal{A}$, it follows that $\mu*a\in I$.

The conditional expectation $\bE$ is $\Gamma$-equivariant in the sense that
\[
   \bE\bigl(\lambda(s) b\,\lambda(s)^*\bigr)
   = \alpha_s\bigl(\bE(b)\bigr),\qquad s\in\Gamma,\ b\in C(X)\rtimes_r\Gamma,
\]
where $\alpha_s(f)=f\circ s^{-1}$ is the induced action on $C(X)$.
Therefore
\[
   \bE(\mu*a)
   = \sum_{s\in H} \mu(s)\,\alpha_s\bigl(\bE(a)\bigr)
   = \mu*\bE(a),
\]
and in particular $\bE(\mu*a)(x) > \delta$ for all $x\in X$.

Now fix $0<\varepsilon<1$. Using \cite[Lemma~4.1]{amrutam2022generalized}
and the assumption that $H$ is $Y$-plump, we can find elements
$g_1,\dots,g_m\in C(Y)$ with $0\le g_j\le 1$ and
$\sum_{j=1}^m g_j^2=1$, together with elements $s_1,\dots,s_m\in H$,
such that
\[
\Biggl\|
   \frac{1}{\delta}\sum_{j=1}^m
      g_j\,\lambda(s_j)\bigl(\mu*a-\bE(\mu*a)\bigr)\lambda(s_j)^* g_j
\Biggr\|
< \varepsilon.
\]
Here we view each $g_j\in C(Y)$ as an element of $C(X)$ via $g_j\circ\pi$,
and hence also as an element of $\mathcal{A}$.

Define
\[
   S := \frac{1}{\delta}\sum_{j=1}^m
         g_j\,\lambda(s_j)\,(\mu*a)\,\lambda(s_j)^* g_j,
   \qquad
   T := \frac{1}{\delta}\sum_{j=1}^m
         g_j\,\lambda(s_j)\,\bE(\mu*a)\,\lambda(s_j)^* g_j.
\]
Then $S\in I$ (because $\mu*a\in I$ and $g_j,\lambda(s_j)\in\mathcal{A}$),
and the above inequality says precisely that
\[
   \|S-T\| < \varepsilon < 1.
\]
Next, we claim that $T$ is invertible in $\mathcal{A}$. For each $x\in X$ we have
\begin{align*}
   T(x)
   &= \frac{1}{\delta}\sum_{j=1}^m g_j(x)^2\,
      \bE(\mu*a)\bigl(s_j^{-1}x\bigr) \\
   &= \frac{1}{\delta}\sum_{j=1}^m g_j(x)^2\,
      (\mu*\bE(a))\bigl(s_j^{-1}x\bigr) \\
   &\ge \frac{1}{\delta}\sum_{j=1}^m g_j(x)^2\,\delta
    = \sum_{j=1}^m g_j(x)^2 = 1,
\end{align*}
where we used $\bE(\mu*a)=\mu*\bE(a)$ and $(\mu*\bE(a))>\delta$ pointwise.
Thus $T$ is a positive element satisfying $T\ge 1$, so it is invertible
and $\|T^{-1}\|\le 1$.

Since $\|S-T\|<1$ and $\|T^{-1}\|\le 1$, we have
\[
   \|T^{-1}(S-T)\| \le \|T^{-1}\|\,\|S-T\| < 1.
\]
Hence $\text{id} + T^{-1}(S-T)$ is invertible, and therefore
$
   S = T\bigl(\text{id} + T^{-1}(S-T)\bigr)
$
is invertible in $C(X)\rtimes_r\Gamma$.

Since $S\in I$ and $S$ is invertible in $\mathcal{A}$, we have
$1=S^{-1}S\in I$, and hence $I=\mathcal{A}$. This shows that
$\mathcal{A}$ is simple.
\end{proof}

Beyond simplicity, one can ask when an intermediate algebra is necessarily a crossed product, a phenomenon known as reflecting~\cite{amrutam2024crossed}. The following theorem provides a criterion for this rigidity in the relative setting, that is, Theorem \ref{thmD}.
\begin{theorem}\label{thm:intmascrossprod} Let $X$ be a $\Gamma$-space, and $\mathcal{B}$ be a unital $\Gamma$-$C^*$-algebra. Consider $C(X)\otimes_{\text{min}}\mathcal{B}$, with the diagonal action.
Let $N = \ker(\Gamma \curvearrowright \mathcal{B})$. Assume that $N$ is $X$-Plump and $N\curvearrowright X$ is minimal. Then for every intermediate C*-algebra $\mathcal{A}$ with
\[
C(X) \rtimes_r \Gamma \subseteq \mathcal{A} \subseteq (C(X) \otimes_{\min} \mathcal{B}) \rtimes_r \Gamma
\]
we have that $E(\mathcal{A}) \subseteq \mathcal{A}$. If, in addition, $\Gamma$ has the approximation property (AP), $\mathcal{A}$ is a crossed product of the form $(C(X) \otimes_{\min} \tilde{\mathcal{B}}) \rtimes_r \Gamma$. Here, $\tilde{\mathcal{B}}$ is a $\Gamma$-invariant $C^*$-subalgebra of $\mathcal{B}$. 
\end{theorem}
\begin{proof}
Let $\mathcal{A}$ be such an intermediate $C^*$-algebra and fix
$a\in\mathcal{A}$. We want to show that $\mathbb{E}(a)\in\mathcal{A}$.

Since $C(X)\otimes_{\min}\mathcal{B}$ is the closed linear span of
finite sums of elementary tensors, we can approximate the coefficient
$\mathbb{E}(a)\in C(X)\otimes_{\min}\mathcal{B}$ by a finite sum:
for any $\varepsilon>0$ there exist $\tilde f_1,\dots,\tilde f_\ell\in C(X)$
and $b_1,\dots,b_\ell\in\mathcal{B}$ such that
\begin{equation}
\label{eq:E(a)}
   \Bigl\|\mathbb{E}(a) - \sum_{j=1}^\ell \tilde f_j\otimes b_j\Bigr\|
   < \frac{\varepsilon}{3}.
\end{equation}
Now fix $x\in X$. Since $N$ is $X$-plump and $N\curvearrowright X$ is
minimal, the argument of \cite[Proposition~3.7]{amrutam2022generalized}
(adapted to the diagonal action on $C(X)\otimes_{\min}\mathcal{B}$ and
using that $N$ acts trivially on $\mathcal{B}$) shows that there are
elements $F_{x,1},\dots,F_{x,m_x}\in C(X)^+$ and
$s_{x,1},\dots,s_{x,m_x}\in N$ with
\[
   \sum_{k=1}^{m_x} F_{x,k}^2 = 1
\]
such that
\begin{equation}
\label{eq:firstineq}
   \Bigl\|
      \sum_{k=1}^{m_x}
         F_{x,k}\,\lambda(s_{x,k})\,a\,\lambda(s_{x,k})^* F_{x,k}
      -(\delta_x\otimes\mathrm{id})(\mathbb{E}(a))
   \Bigr\| < \frac{\varepsilon}{3}.
\end{equation}
Here $\delta_x:C(X)\to\mathbb C$ is evaluation at $x$, and we regard
$(\delta_x\otimes\mathrm{id})(\mathbb{E}(a))$ as an element
of $C(X)\otimes_{\min}\mathcal{B}$ (and hence of the crossed product) via
$1_X\otimes\mathcal{B}$.
Moreover, we also have that
\begin{equation}
\label{eq:sliced}
\Bigl\|\delta_x\otimes\text{id}(\mathbb{E}(a)) - \sum_{j=1}^\ell \tilde f_j(x)\otimes b_j\Bigr\|
   < \frac{\varepsilon}{3}    
\end{equation}
For this fixed $x$, by continuity of each $\tilde f_j$ we can choose an
open neighborhood $U_x\ni x$ such that for all $y\in U_x$ and all
$j=1,\dots,\ell$,
\[
   \|\tilde f_j(y)-\tilde f_j(x)\| < \frac{\varepsilon}{3}.
\]
Since $X$ is compact, there exist finitely many points
$x_1,\dots,x_p\in X$ such that $X=\bigcup_{r=1}^p U_{x_r}$. Let
$\{\tilde F_r\}_{r=1}^p\subset C(X)^+$ be a partition of unity subordinate
to this cover, i.e.\ $0\le \tilde F_r\le 1$, $\supp(\tilde F_r)\subset U_{x_r}$
and $\sum_{r=1}^p \tilde F_r^2 = 1$.

For each $r$ we have functions $F_{x_r,k}\in C(X)^+$ and elements
$s_{x_r,k}\in N$ as above; set
\[
   G_{r,k} := \tilde F_r\,F_{x_r,k},\qquad
   t_{r,k} := s_{x_r,k}.
\]
Then
\[
   \sum_{r,k} G_{r,k}^2
   = \sum_{r=1}^p \tilde F_r^2 \sum_{k=1}^{m_{x_r}} F_{x_r,k}^2
   = \sum_{r=1}^p \tilde F_r^2 = 1.
\]
Define
\[
   S := \sum_{r,k} G_{r,k}\,\lambda(t_{r,k})\,a\,
              \lambda(t_{r,k})^* G_{r,k}.
\]
Since $a\in\mathcal{A}$ and each $G_{r,k},\lambda(t_{r,k})\in
C(X)\rtimes_r\Gamma\subseteq\mathcal{A}$, we have $S\in\mathcal{A}$. By the choice of $F_{x_r,k}$ and $s_{x_r,k}$, together with the fact
that $t_{r,k}\in N$ acts trivially on $\mathcal{B}$, using equation~\eqref{eq:firstineq}, one checks (as in
the proof of \cite[Theorem~1.2]{amrutam2022generalized}) that
\begin{equation}
\label{eq:partitionofunityplmus}
   \bigl\|S - \sum_{r=1}^p \tilde F_r
                    (\delta_{x_r}\otimes\mathrm{id})(\mathbb{E}(a))
                    \tilde F_r\bigr\|
   < \frac{\epsilon}{3}.
\end{equation}
Moreover, using the choice of the neighborhoods $U_{x_r}$ and partition of
unity along with equation~\eqref{eq:sliced}, we see that
\begin{equation}
\label{eq:secondineq}
 \Bigl\|
      \sum_{r=1}^p \tilde F_r
         (\delta_{x_r}\otimes\mathrm{id})(\mathbb{E}(a)) \tilde F_r
      - \sum_{j=1}^\ell \tilde f_j\otimes b_j
   \Bigr\|
   < \frac{\epsilon}{3}.   \end{equation}
Combining equations~\eqref{eq:E(a)},~\eqref{eq:partitionofunityplmus} and~\eqref{eq:secondineq},
we obtain that
\begin{align*}
  \|S-\mathbb{E}(a)\|
\le&
 \Biggl\|
  S
   - \sum_{r=1}^p \tilde F_r
      (\delta_{x_r}\otimes\mathrm{id})(\mathbb{E}(a))\,
      \tilde F_r
 \Biggr\|\\
 & +
 \Biggr\|  \sum_{r=1}^p \tilde F_r
      (\delta_{x_r}\otimes\mathrm{id})(\mathbb{E}(a))\,
      \tilde F_r-
\sum_{j=1}^\ell \tilde f_j\otimes b_j
 \Biggr\|
 + \Biggl\|\sum_{j=1}^\ell \tilde f_j\otimes b_j - \mathbb{E}(a)\Biggr\| \\
& <\varepsilon.
\end{align*}
Since $\varepsilon>0$ was arbitrary and $\mathcal{A}$ is norm closed, it
follows that $\mathbb{E}(a)\in \mathcal{A}$.

Thus $\mathbb{E}(\mathcal{A})\subseteq\mathcal{A}$. If, in addition,
$\Gamma$ has  AP, then by
\cite[Proposition~3.4]{suzuki2017group} any intermediate
$C^*$-subalgebra between $C(X)\rtimes_r\Gamma$ and
$(C(X)\otimes_{\min}\mathcal{B})\rtimes_r\Gamma$ that is invariant under
$\mathbb{E}$ must be a crossed product. Hence $\mathcal{A}$ is a crossed
product of the form $\tilde{A}\rtimes_r\Gamma$, where $C(X)\subset\tilde{\mathcal{A}}\subset C(X)\otimes_{\text{min}}\mathcal{B}$. Since $N\curvearrowright X$ is minimal, $C(X)$ is $N$-simple. Moreover, $N\curvearrowright\mathcal{B}$ is trivial. We can now appeal to \cite[Proposition~1.2]{amrutam2025splitting} to obtain that $\tilde{\mathcal{A}}=C(X)\otimes_{\text{min}}\tilde{\mathcal{B}}$. The proof is complete.  
\end{proof}
Unlike the examples provided in \cite{amrutam2021intermediate, ursu2022relative}, which typically require the ambient group $\Gamma$ to be $C^*$-simple, our construction allows the groups to be non-$C^*$-simple. Consequently, we obtain rigidity results for a broad spectrum of groups, including those with non-trivial amenable radicals, providing a unified approach that encompasses both standard product groups and more complex extensions, such as wreath products.

To prepare for our examples, let us briefly recall the following lemma, 
which can in fact be obtained from Lemma~3.2  together with the argument in the proof of Theorem~3.1 in \cite{breuillard2017c}. However, for the reader's convenience, we provide a self-contained proof.
\begin{lemma}
\label{lem:topfree-lift}
Let $X$ and $Z$ be 
$\Gamma$-spaces.
Suppose that $\rho:Z\to X$ is a  $\Gamma$-factor map, that
$\Gamma\curvearrowright Z$ is minimal, and that $\Gamma\curvearrowright X$
is topologically free.
Then $\Gamma\curvearrowright Z$ is topologically free.
\end{lemma}
\begin{proof}
Fix $s\in\Gamma\setminus\{e\}$ and consider the fixed-point set
\[
   \Fix_Z(s) := \{z\in Z : sz=z\}.
\]
Suppose for a contradiction that $\Fix_Z(s)$ has non-empty interior.
Then there exists a non-empty open set $U\subset \Fix_Z(s)$.

Since $\Gamma\curvearrowright Z$ is minimal, it follows from \cite[Theorem 1.15]{Auslanderbook} (see also \cite[Lemma 3.2]{breuillard2017c} or \cite[Lemma 2.8]{MR4348698}) that $\rho(U)$ has non-empty interior. Note that $\rho(U)\subset  \Fix_X(s) $, a contradiction. Thus,  $\Gamma\curvearrowright Z$ is topologically free.
\end{proof}
We now provide a recipe for producing $X$-plump subgroups
without assuming that the ambient group $\Gamma$ is $C^*$-simple.
\begin{example}\label{exa:diagonal-plump}
Let $\Gamma$ be a discrete group and let $N\lhd\Gamma$ be a non-trivial
normal subgroup.
Assume that $N$ is non-amenable (so that the Furstenberg boundary
$\partial_F N$ is non-trivial).
Let $X$ be a $\Gamma$-space such that
\begin{enumerate}[(i)]
  \item the restricted action $N\curvearrowright X$ is minimal;
  \item the action $\Gamma\curvearrowright X$ is topologically free.
\end{enumerate}

Consider the product space
\[
   Y := X\times \partial_F N
\]
with the diagonal $N$-action
\[
   n\cdot(x,\xi) := (nx, n\xi),\qquad n\in N,\,x\in X,\,\xi\in\partial_F N.
\]
Since the action $N\curvearrowright\partial_F N$ extends to a $\Gamma$-action~\cite[Lemma~5.2]{breuillard2017c},
we may view $Y$ as a $\Gamma$-space via the diagonal $\Gamma$-action. By Example~\ref{Ex:strongproximal}, there exists a unique closed
$\Gamma$-invariant subset $\widetilde X\subset Y$ such that
\begin{itemize}
  \item $N\curvearrowright \widetilde X$ is minimal, and
  \item the restriction of the first-coordinate projection
        \[
           \pi : \widetilde X\to X
        \]
        is onto and is a generalized boundary for the pair $(N,X)$;
        in particular, $N\curvearrowright\widetilde X$ is a strongly
        proximal extension over $X$.
\end{itemize}
By Lemma~\ref{lem:topfree-lift}, the action $\Gamma\curvearrowright \widetilde X$ is topologically free. Consequently, we get that $\pi:\widetilde X\to X$ is a generalized boundary for $(N,X)$
such that
\begin{itemize}
  \item $N\curvearrowright \widetilde X$ is a strongly proximal extension
        over $X$ and is minimal;
  \item $\Gamma\curvearrowright \widetilde X$ is topologically free.
\end{itemize}
Hence, by Theorem~\ref{thm:condgenplump}, the subgroup $N$ is $X$-plump.
\end{example}
Putting together Theorem~\ref{thm:intmascrossprod} along with Example~\ref{exa:diagonal-plump}, we obtain the following.
\begin{corollary}
Let $\Gamma$ be a non-amenable group (not necessarily $C^*$-simple) acting
non-faithfully on a unital $C^*$-algebra $\mathcal{B}$, and let
$ N := \ker(\Gamma\curvearrowright \mathcal{B})$.
Let $X$ be a  $\Gamma$-space such that $X$ is $N$-minimal and
$\Gamma\curvearrowright X$ is topologically free.
Consider the $\Gamma$-$C^*$-algebra $C(X)\otimes_{\min}\mathcal{B}$ with
the diagonal action.
Then for every intermediate $C^*$-algebra $\mathcal{A}$ with
\[
   C(X)\rtimes_r\Gamma
   \;\subseteq\;
   \mathcal{A}
   \;\subseteq\;
   (C(X)\otimes_{\min}\mathcal{B})\rtimes_r\Gamma
\]
we have $E(\mathcal{A})\subseteq\mathcal{A}$.
If, in addition, $\Gamma$ has the approximation property \emph{(AP)}, then
$\mathcal{A}$ is a crossed product of the form $(C(X) \otimes_{\min} \tilde{\mathcal{B}}) \rtimes_r \Gamma$.
\end{corollary}
\begin{proof}
By Example~\ref{exa:diagonal-plump}, the subgroup $N$ is $X$-plump. Since $N\curvearrowright X$ is minimal by assumption, the result follows immediately from Theorem~\ref{thm:intmascrossprod}.    
\end{proof}

The following is from a combination of Corollary~\ref{cor:connected} and Theorem~\ref{thm:intmascrossprod}.    
\begin{corollary}Let $\Gamma$ be a countable group acting non-faithfully on a unital
$C^*$-algebra $\mathcal{B}$, and set
$
   N := \ker(\Gamma \curvearrowright \mathcal{B}).
$
Assume that 
$C_{\Gamma}(N)=\{e\}$.
Let $X$ be a connected  $\Gamma$-space such that the restricted
action $N \curvearrowright X$ is minimal, and assume that the reduced
crossed product $C(X)\rtimes_r N$ is simple.
Consider the $\Gamma$-$C^*$-algebra $C(X)\otimes_{\min}\mathcal{B}$ with
the diagonal action.
Then for every intermediate $C^*$-algebra $\mathcal{A}$ with
\[
   C(X) \rtimes_r \Gamma
   \;\subseteq\;
   \mathcal{A}
   \;\subseteq\;
   (C(X) \otimes_{\min} \mathcal{B}) \rtimes_r \Gamma
\]
we have $E(\mathcal{A}) \subseteq \mathcal{A}$.
If, in addition, $\Gamma$ has AP, then
$\mathcal{A}$ is a crossed product of the form $(C(X) \otimes_{\min} \tilde{\mathcal{B}}) \rtimes_r \Gamma$. 
\end{corollary}
\begin{proof}
Since $N$ is the kernel of the action on $\mathcal{B}$, it is a normal subgroup of $\Gamma$, and hence $N$ is commensurated in $\Gamma$ (i.e., $N \le_c \Gamma$). By Corollary~\ref{cor:connected}, the assumptions that $X$ is connected, $N\curvearrowright X$ is minimal, $C(X) \rtimes_r N$ is simple, and $C_{\Gamma}(N)=\{e\}$ imply that $N$ is $X$-plump. The conclusion then follows directly from Theorem~\ref{thm:intmascrossprod}.    
\end{proof}
\bibliography{name}
\bibliographystyle{amsalpha}

\end{document}